\documentclass{agtart_a}
\pdfoutput=1

\usepackage{amscd}
\usepackage[matrix,arrow,curve,frame]{xy}    % XY-pic diagram pac

\title[Continuous functors as a model for the equivariant stable homotopy category]{Continuous functors as a model\\for the equivariant stable homotopy category}

\author{Andrew J\,Blumberg}
\givenname{Andrew}
\surname{Blumberg}
\address{Department of Mathematics\\
Stanford University\\\newline
450 Serra Mall\\
Stanford, California 94305\\USA} 
\email{blumberg@math.stanford.edu}
\urladdr{}

\volumenumber{6}
\issuenumber{}
\publicationyear{2006}
\papernumber{79}
\startpage{2257}
\endpage{2295}

\doi{}
\MR{}
\Zbl{}

\keyword{equivariant infinite loop space theory}
\keyword{excisive functors}
\keyword{Atiyah duality}
\subject{primary}{msc2000}{55P42}
\subject{secondary}{msc2000}{55P47}
\subject{secondary}{msc2000}{55P91}

\received{15 June 2005}
\revised{8 November 2006}
\accepted{9 November 2006}
\published{8 December 2006}
\publishedonline{8 December 2006}
\proposed{}
\seconded{}
\corresponding{}
\editor{Liz}
\version{}

\arxivreference{math.AT/0505512}

\let\xysavmatrix\xymatrix
\def\xymatrix{\disablesubscriptcorrection\xysavmatrix}
\AtBeginDocument{\let\bar\wbar\let\tilde\wtilde\let\hat\what}
\newcommand{\sW}{\scr{W}\!}
\newcommand{\sWuG}{\scr{W}_G}
\newcommand{\sWuH}{\scr{W}_H}
\newcommand{\sT}{\scr{T}\!}
\newcommand{\sF}{\scr{F}\!}
\newcommand{\sI}{\scr{I}\!}
\newcommand{\sS}{\scr{S}\!}

\DeclareFontFamily{OMS}{rsfs}{\skewchar\font'60}
\DeclareFontShape{OMS}{rsfs}{m}{n}{<-5>rsfs5 <5-7>rsfs7 <7->rsfs10 }{}
\DeclareSymbolFont{rsfs}{OMS}{rsfs}{m}{n}
\DeclareSymbolFontAlphabet{\scr}{rsfs}

\newcommand{\sC}{\scr{C}}

\newcommand{\sP}{\scr{P}}

\newcommand{\bP}{\mathbb{P}}

\newcommand{\bU}{\mathbb{U}}

\newcommand{\htp}{\simeq}    % homotopy symbol
\newcommand{\thp}{\ltimes}   % twisted half-smash product
\newcommand{\sma}{\wedge}    % smash product 
\newcommand{\monoto}{\lhook\joinrel\relbar\joinrel\rightarrow}

\makeatletter
\def\cnewtheorem#1[#2]#3{\newtheorem{#1}{#3}[section]
\expandafter\let\csname c@#1\endcsname\c@thm}
\makeatother

\newtheorem{thm}{Theorem}[section]

\cnewtheorem{cor}[thm]{Corollary}
\cnewtheorem{corollary}[thm]{Corollary}
\cnewtheorem{prop}[thm]{Proposition}
\cnewtheorem{proposition}[thm]{Proposition}
\cnewtheorem{lem}[thm]{Lemma}
\cnewtheorem{lemma}[thm]{Lemma}
\cnewtheorem{conj}[thm]{Conjecture}
\cnewtheorem{prob}[thm]{Problem}
\cnewtheorem{exer}[thm]{Exercise}
\cnewtheorem{hyp}[thm]{Hypothesis}
\makeautorefname{hyp}{Hypothesis}
\cnewtheorem{convention}[thm]{Convention}

\theoremstyle{definition}
\cnewtheorem{defn}[thm]{Definition}
\cnewtheorem{definition}[thm]{Definition}
\cnewtheorem{defns}[thm]{Definitions}
\cnewtheorem{exmp}[thm]{Example}
\cnewtheorem{exmps}[thm]{Examples}
\cnewtheorem{var}[thm]{Variant}
\cnewtheorem{vars}[thm]{Variants}
\cnewtheorem{con}[thm]{Construction}
\cnewtheorem{notn}[thm]{Notation}
\cnewtheorem{notns}[thm]{Notations}
\cnewtheorem{conv}[thm]{Convention}

\theoremstyle{remark}
\cnewtheorem{rem}[thm]{Remark}
\cnewtheorem{remark}[thm]{Remark}
\cnewtheorem{rems}[thm]{Remarks}
\cnewtheorem{warn}[thm]{Warning}
\cnewtheorem{sch}[thm]{Scholium}
\cnewtheorem{expl}[thm]{Explanations}
\cnewtheorem{ouch}[thm]{Counterexample}

\newcommand{\colim}{\ensuremath{\operatorname{colim}}}

\newcommand{\id}{\ensuremath{\operatorname{id}}}

\newcommand{\Map}{\ensuremath{\operatorname{Map}}}

\begin{document}

\begin{asciiabstract}
It is a classical observation that a based continuous functor X from
the category of finite CW--complexes to the category of based spaces
that takes homotopy pushouts to homotopy pullbacks ``represents'' a
homology theory---the collection of spaces {X(S^n)} obtained by
evaluating X on spheres yields an Omega-prespectrum.  Such functors
are sometimes referred to as linear or excisive.  The main theorem of
this paper provides an equivariant analogue of this result. We show
that a based continuous functor from finite G-CW--complexes to based
G-spaces represents a genuine equivariant homology theory if and only
if it takes G-homotopy pushouts to G-homotopy pullbacks and satisfies
an additional condition requiring compatibility with Atiyah duality
for orbit spaces G/H.

Our motivation for this work is the development of a recognition
principle for equivariant infinite loop spaces.  In order to make the
connection to infinite loop space theory precise, we reinterpret the
main theorem as providing a fibrancy condition in an appropriate model
category of spectra.  Specifically, we situate this result in the
context of the study of equivariant diagram spectra indexed on the
category W_G of based G-spaces homeomorphic to finite G-CW--complexes
for a compact Lie group G.  Using the machinery of
Mandell--May--Schwede--Shipley, we show that there is a stable model
structure on this category of diagram spectra which admits a monoidal
Quillen equivalence to the category of orthogonal G-spectra.  We
construct a second ``absolute'' stable model structure which is
Quillen equivalent to the stable model structure.  There is a
model-theoretic identification of the fibrant continuous functors in
the absolute stable model structure as functors Z such that for A in
W_G the collection {Z(A smash S^W)} forms an Omega-G-prespectrum as W
varies over the universe U.  Thus, our main result provides a concrete
identification of the fibrant objects in the absolute stable model
structure.

This description of fibrant objects in the absolute stable model
structure makes it clear that in the equivariant setting we cannot
hope for a comparison between the category of equivariant continuous
functors and equivariant Gamma-spaces, except when G is finite. We
provide an explicit analysis of the failure of the category of
equivariant Gamma-spaces to model connective G-spectra, even for G =
S^1.
\end{asciiabstract}

\begin{htmlabstract}
<p class="noindent">
It is a classical observation that a based continuous functor X from
the category of finite CW&ndash;complexes to the category of based spaces
that takes homotopy pushouts to homotopy pullbacks &ldquo;represents&rdquo; a
homology theory&ndash;-the collection of spaces {X(S<sup>n</sup>)} obtained by
evaluating X on spheres yields an &Omega;&ndash;prespectrum.  Such
functors are sometimes referred to as linear or excisive.  The main
theorem of this paper provides an equivariant analogue of this
result. We show that a based continuous functor from finite
G&ndash;CW&ndash;complexes to based G&ndash;spaces represents a genuine
equivariant homology theory if and only if it takes G&ndash;homotopy
pushouts to G&ndash;homotopy pullbacks and satisfies an additional
condition requiring compatibility with Atiyah duality for orbit spaces
G/H.
</p>
<p class="noindent">
Our motivation for this work is the development of a recognition
principle for equivariant infinite loop spaces.  In order to make the
connection to infinite loop space theory precise, we reinterpret the
main theorem as providing a fibrancy condition in an appropriate model
category of spectra.  Specifically, we situate this result in the
context of the study of equivariant diagram spectra indexed on the
category W<sub>G</sub> of based G&ndash;spaces homeomorphic to finite
G&ndash;CW&ndash;complexes for a compact Lie group G.  Using the machinery
of Mandell&ndash;May&ndash;Schwede&ndash;Shipley, we show that there is a stable
model structure on this category of diagram spectra which admits a
monoidal Quillen equivalence to the category of orthogonal
G&ndash;spectra.  We construct a second &ldquo;absolute&rdquo; stable model
structure which is Quillen equivalent to the stable model structure.
There is a model-theoretic identification of the fibrant continuous
functors in the absolute stable model structure as functors Z such
that for A &isin; W<sub>G</sub> the collection {Z(A &and; S<sup>W</sup>)} forms an
&Omega;&ndash;G&ndash;prespectrum as W varies over the universe U.  Thus,
our main result provides a concrete identification of the fibrant
objects in the absolute stable model structure.
</p>
<p class="noindent">
This description of fibrant objects in the absolute stable model
structure makes it clear that in the equivariant setting we cannot
hope for a comparison between the category of equivariant continuous
functors and equivariant &Gamma;&ndash;spaces, except when G is finite.
We provide an explicit analysis of the failure of the category of
equivariant &Gamma;&ndash;spaces to model connective G&ndash;spectra, even
for G = S<sup>1</sup>.
</p>
\end{htmlabstract}

\begin{webabstract}
It is a classical observation that a based continuous functor $X$ from
the category of finite CW--complexes to the category of based spaces
that takes homotopy pushouts to homotopy pullbacks ``represents'' a
homology theory---the collection of spaces $\{X(S^n)\}$ obtained by
evaluating $X$ on spheres yields an $\Omega$--prespectrum.  Such
functors are sometimes referred to as linear or excisive.  The main
theorem of this paper provides an equivariant analogue of this
result. We show that a based continuous functor from finite
$G$--CW--complexes to based $G$--spaces represents a genuine
equivariant homology theory if and only if it takes $G$--homotopy
pushouts to $G$--homotopy pullbacks and satisfies an additional
condition requiring compatibility with Atiyah duality for orbit spaces
$G/H$.

Our motivation for this work is the development of a recognition
principle for equivariant infinite loop spaces.  In order to make the
connection to infinite loop space theory precise, we reinterpret the
main theorem as providing a fibrancy condition in an appropriate model
category of spectra.  Specifically, we situate this result in the
context of the study of equivariant diagram spectra indexed on the
category $\mathcal{W}_G$ of based $G$--spaces homeomorphic to finite
$G$--CW--complexes for a compact Lie group $G$.  Using the machinery
of Mandell--May--Schwede--Shipley, we show that there is a stable
model structure on this category of diagram spectra which admits a
monoidal Quillen equivalence to the category of orthogonal
$G$--spectra.  We construct a second ``absolute'' stable model
structure which is Quillen equivalent to the stable model structure.
There is a model-theoretic identification of the fibrant continuous
functors in the absolute stable model structure as functors $Z$ such
that for $A \in \mathcal{W}_G$ the collection $\{Z(A \wedge S^W)\}$ forms an
$\Omega$--$G$--prespectrum as $W$ varies over the universe $U$.  Thus,
our main result provides a concrete identification of the fibrant
objects in the absolute stable model structure.

This description of fibrant objects in the absolute stable model
structure makes it clear that in the equivariant setting we cannot
hope for a comparison between the category of equivariant continuous
functors and equivariant $\Gamma$--spaces, except when $G$ is finite.
We provide an explicit analysis of the failure of the category of
equivariant $\Gamma$--spaces to model connective $G$--spectra, even
for $G = S^1$.
\end{webabstract}

\begin{abstract}
It is a classical observation that a based continuous functor $X$ from
the category of finite CW--complexes to the category of based spaces
that takes homotopy pushouts to homotopy pullbacks ``represents'' a
homology theory---the collection of spaces $\{X(S^n)\}$ obtained by
evaluating $X$ on spheres yields an $\Omega$--prespectrum.  Such
functors are sometimes referred to as linear or excisive.  The main
theorem of this paper provides an equivariant analogue of this
result. We show that a based continuous functor from finite
$G$--CW--complexes to based $G$--spaces represents a genuine
equivariant homology theory if and only if it takes $G$--homotopy
pushouts to $G$--homotopy pullbacks and satisfies an additional
condition requiring compatibility with Atiyah duality for orbit spaces
$G/H$.

Our motivation for this work is the development of a recognition
principle for equivariant infinite loop spaces.  In order to make the
connection to infinite loop space theory precise, we reinterpret the
main theorem as providing a fibrancy condition in an appropriate model
category of spectra.  Specifically, we situate this result in the
context of the study of equivariant diagram spectra indexed on the
category $\mathscr{W}_G$ of based $G$--spaces homeomorphic to finite
$G$--CW--complexes for a compact Lie group $G$.  Using the machinery
of Mandell--May--Schwede--Shipley, we show that there is a stable
model structure on this category of diagram spectra which admits a
monoidal Quillen equivalence to the category of orthogonal
$G$--spectra.  We construct a second ``absolute'' stable model
structure which is Quillen equivalent to the stable model structure.
There is a model-theoretic identification of the fibrant continuous
functors in the absolute stable model structure as functors $Z$ such
that for $A \in \mathscr{W}_G$ the collection $\{Z(A \wedge S^W)\}$ forms an
$\Omega$--$G$--prespectrum as $W$ varies over the universe $U$.  Thus,
our main result provides a concrete identification of the fibrant
objects in the absolute stable model structure.

This description of fibrant objects in the absolute stable model
structure makes it clear that in the equivariant setting we cannot
hope for a comparison between the category of equivariant continuous
functors and equivariant $\Gamma$--spaces, except when $G$ is finite.
We provide an explicit analysis of the failure of the category of
equivariant $\Gamma$--spaces to model connective $G$--spectra, even
for $G = S^1$.
\end{abstract}

\maketitle

\section{Introduction}

One of the striking successes in the development of stable homotopy theory was the 
characterization of infinite loop spaces, spaces that arise as the zero space of a spectrum.  
Following Boardman and Vogt \cite{boardman-vogt}, the approaches of May via 
$E_\infty$--operads \cite{may-geom} and Segal via $\Gamma$--spaces \cite{segal} provided 
characterizations of space-level data guaranteeing that a space possessed arbitrary 
deloopings.  However, while in general the development of stable homotopy theory in the 
equivariant setting has been successful as in Lewis, May and Steinberger \cite{lms}, the 
area of equivariant infinite loop space theory has remained mysterious.

In the nonequivariant setting, an infinite loop space is a space equipped with a 
multiplication which is commutative and associative up to all higher homotopies.  The 
recognition principles explicitly encode this information---both $E_\infty$--operad actions 
and $\Gamma$--space structures are evidently devices for packaging up the coherent homotopies 
describing such a multiplication.  Unfortunately, in the equivariant setting the structure 
carried by an infinite loop space is harder to understand.  The additional complexity arises 
from the representation theory of the group $G$.

For certain applications, one can work with equivariant spectra which consist of a sequence 
of $G$--spaces $X_n$ and equivariant structure maps $S^1 \sma X_n \rightarrow X_{n+1}$ which 
induce $G$--equivalences $X_n \rightarrow \Omega X^{n+1}$.  However, in order to have an 
equivariant version of Spanier--Whitehead duality, one has no choice but to work with spectra 
indexed on the collection of all finite-dimensional real representations of the group $G$ and 
equipped with structure maps for suspensions along such representations.

For a finite-dimensional real representation $V$, let $S^V$ denote the one-point 
compactification.  A ``genuine'' equivariant spectrum $X$ has compatible structure maps 
$S^{W-V} \sma X(V) \rightarrow X(W)$ which induce adjoint equivalences $X(V) \rightarrow 
\Omega^{W-V} X(W)$ for all pairs of finite-dimensional representations $V \subset W$, where 
$W-V$ is the orthogonal complement of $V$ in $W$.  The zero space of such a spectrum carries 
a tremendous amount of structure inherited from the group $G$---it is a ``$`V$--fold'' loop 
space for all finite-dimensional representations $V$.

Nonetheless, when $G$ is a finite group, straightforward generalizations of the 
nonequivariant recognition principles continue to apply.  In the operadic setting, there is a 
notion of a $G$--operad and one can demonstrate that spaces that admit the action of an 
$E_\infty$--$G$--operad admit deloopings by arbitrary representation spheres; see Costenoble 
and Waner \cite{costenoble-waner}.  The equivariant argument follows the nonequivariant 
argument closely, and in particular depends on equivariant versions of the ``approximation 
theorem'', which describes models for the free loop spaces $\Omega^V S^V X$.  Similarly, one 
can define an equivariant version of a $\Gamma$--space, and it can be shown that a ``very 
special'' equivariant $\Gamma$--space gives rise to a genuine $G$--spectrum; see Shimakawa 
\cite{shimakawa} and Segal \cite{equi-segal}.

Unfortunately, serious difficulties arise in trying to generalize to the situation in which 
$G$ is an infinite compact Lie group.  There is a simple counterexample due to Segal which 
demonstrates that the equivariant version of the approximation theorem fails, even for $G = 
S^1$ \cite{equi-segal}.  There are problems in trying to generalize the $\Gamma$--space 
approach as well.  In the appendix, we will recall the counterexample of Segal and also 
analyze the failure of equivariant $\Gamma$--spaces to be a model for the equivariant stable 
homotopy category for any infinite compact Lie group $G$.

The intent of this paper is to begin to characterize the kind of structure that arises on the 
zero spaces of genuine equivariant spectra by studying the closely related category of 
equivariant continuous functors.  This approach is suggested by the strategy of Segal's 
analysis of $\Gamma$--spaces.  Recall that Segal proved that $\Gamma$--spaces describe 
infinite loop spaces by showing that a $\Gamma$--space gives rise via prolongation to a 
continuous functor from finite CW--complexes to based spaces which takes (most) homotopy 
pushouts to homotopy pullbacks \cite{segal}.  Let $\sT$ denote the category of compactly 
generated based spaces and $\sW$ denotes the full subcategory of $\sT$ consisting of spaces 
homeomorphic to finite CW--complexes.  A continuous functor $X \co \sW \rightarrow \sT$ is 
one for which the induced map of hom spaces
\[
\xymatrix{ \relax\sW(A, B) \ar[r]^-= & \relax \sT(A, B) \ar[r]^-X &
\relax\sT(X(A), X(B)) }
\]
is continuous.  A continuous functor $X$ between pointed categories is based if $X(*) = *$.  
We will assume herein that all continuous functors are based.  It is reasonable to think of a 
continuous functor as a kind of spectrum because continuity implies the existence of a 
``structure map'' $X(A) \sma B \rightarrow X(A \sma B)$.  Indeed, in the nonequivariant 
setting Lydakis showed that a simplicial version of this category was a convenient symmetric 
monoidal category of spectra \cite{lydakis}.  When a continuous functor $X$ takes homotopy 
pushouts to homotopy pullbacks, then the prespectrum $\{X(S^n)\}$ is an 
$\Omega$--prespectrum.  Applying $X$ to the diagram 
\[
\xymatrix{
A \ar[r]\ar[d] & \relax* \ar[d] \\
\relax* \ar[r] & S^1 \sma A \\
}
\]
induces a weak equivalence between $X(A)$ and $\Omega X(\Sigma A)$.
Therefore, we are led to seek the equivariant generalization of this
property.  In order to make this question precise, we need to specify
exactly what we mean by an equivariant prespectrum.

Fix a universe $U$, by which we mean an infinite (countable) dimensional real vector space 
$U$ on which $G$ acts by isometries and which is the direct sum of its finite-dimensional 
$G$--invariant subspaces.  We will assume that $U$ contains a trivial representation and each 
of its finite-dimensional subrepresentations infinitely often.  The universe $U$ is complete 
if it contains all irreducible representations of $G$.  

\begin{remark}
In the body of the paper, for simplicity we will assume that the universes we work with are 
complete.  The results of the paper, suitably modified, remain true for incomplete universes.  
The specific modifications necessary and related subtleties which arise in the context of 
incomplete universes are discussed in \fullref{s:incomplete}.
\end{remark}

A $G$--prespectrum is a collection of spaces $X(V)$ for finite-dimensional $V \in U$ equipped 
with compatible structure maps $S^{W-V} \sma X(V) \rightarrow X(W)$.  An 
$\Omega$--$G$--prespectrum is a $G$--prespectrum $X$ with adjoint structure maps $X(V) 
\rightarrow \Omega^{W-V} X(W)$ which are $G$--equivalences.

Thus, for a compact Lie group $G$ we wish to know when a continuous functor $Z$ from finite 
$G$--CW--complexes to $G$--spaces yields a collection of spaces $\{Z(V)\}$ for $Z \in U$ 
which specifies a genuine equivariant $\Omega$--prespectrum indexed on the universe $U$. It 
is not sufficient for such a $Z$ to take $G$--homotopy pullbacks to $G$--homotopy pushouts, 
as the structure maps for nontrivial representations $V$ cannot be constructed in the above 
fashion.

Denote the category of based $G$--spaces that are homeomorphic to finite $G$--CW--complexes 
by $\sWuG$.  Note that we include all maps as morphisms, not just the equivariant maps.  As a 
consequence, $\sWuG$ is enriched over based $G$--spaces with the action on the morphism 
spaces given by conjugation.  In analogy with the nonequivariant terminology, we say a 
functor from $\sWuG$ to based $G$--spaces as continuous if the induced map of enriched $\hom$ 
$G$--spaces is a continuous map of $G$--spaces.  We will refer to a based continuous functor 
from $\sWuG$ to based $G$--spaces as a $\sWuG$--space.

The main result of this paper is the following theorem, which implies that the additional 
sufficient condition for a $\sWuG$--space $Z$ to represent an $\Omega$--$G$--prespectrum is a 
kind of compatibility with Spanier--Whitehead duality for the orbit spectra $\Sigma^\infty 
G/H_+$. Let $G/H$ be embedded in a real $G$--representation $V$, with normal bundle $\nu$.  
Denote by $T\nu$ the Thom space of $\nu$.  Recall that there is a stable equivalence between 
$\Sigma^\infty G/H_+$ and $\Sigma^\infty_V T\nu$ ($S^{-V} \sma T\nu$) as a consequence of 
Spanier--Whitehead duality.  This stable duality is exhibited on the space level by a 
$V$--duality map $G/H_+ \sma T\nu \rightarrow S^V$.

\begin{thm}\label{t:char}
Let $G$ be a compact Lie group, $U$ a complete universe of $G$--represent\-a\-tions, and $Z$ 
a $\sWuG$--space.  Let $A$ be any finite $G$--CW--complex.  Then the collection $\{Z(A \sma 
S^V) \, | \, V \subset U\}$ is an $\Omega$--$G$--prespectrum if and only if $(1)$ and either 
of the equivalent conditions $(2)$ or $\mathrm{(2^{\prime})}$ hold.

\begin{enumerate}
\item{$Z$ takes $G$--homotopy pushout squares to $G$--homotopy pullback squares.}
\item{Let $G/H$ be an orbit space embedded in a $G$--representation $V
    \subset U$, with normal bundle $\nu$.  Let $T\nu$ denote the Thom
    space of $\nu$.  For any $X \in \sWuG$, a certain
    map \[Z(T\nu \sma X) \rightarrow \Map_0(G/H_+, Z(S^V \sma X))\] is
    an equivalence.  Here $\Map_0$ denotes the $G$--space of
    nonequivariant based maps, with $G$ acting by conjugation.}  
\item[$\mathrm{(2^{\prime})}$]{Let $G/H$ be an orbit space embedded in a
    $G$--representation $V \subset U$.  Let $L$ denote the tangent
    $H$--representation at the identity coset.  For any $X \in \sWuH$,
    a certain map \[Z(G_+ \sma_H X) \rightarrow
    \Map_H(G_+, Z(S^L \sma X))\] is an equivalence.} 
\end{enumerate}

The maps in conditions $(2)$ and $(2^{\prime})$ are induced from the $V$--duality $T\nu \sma 
G/H_+ \rightarrow S^V$ and will be described precisely in \fullref{d:xi}.  The second version 
of the orbit condition arises from the generalized Wirthmuller isomorphism \cite{may,lms}.
\end{thm}

A based continuous functor $Z$ which satisfies either of the
equivalent conditions above will be said to be ``equivariantly linear
with respect to $U$''.  This terminology is motivated by the fact that
our notion of equivariantly linear provides a precise generalization
of the linearity conditions of Mandell, May, Schwede and Shipley \cite{mmss}.  One might also describe
such a functor $Z$ as ``genuinely'' equivariantly excisive, in line
with the language of Goodwillie's calculus of functors.  It is
interesting to wonder if this notion is relevant to possible
equivariant generalizations of the calculus of functors.

Before we move on, it is worth attempting to provide some intuition about why this 
compatibility with Atiyah duality for orbit spectra is a reasonable condition to expect, 
beyond what is provided by the details of the proof of the theorem.  A related problem to the 
one we consider herein is to determine when a $\Z$--graded cohomology theory on $G$--spaces 
with coefficients in a coefficient system $M$ extends to an $RO(G)$--graded cohomology 
theory.  Such an extension exists if and only if $M$ extends to a Mackey functor 
\cite[IX.5.2]{may-alaska}. It is illuminating to recall the data that is required for such an 
extension.  Essentially, the key observation is that the stable transfer maps $\tau (G/H) \co 
S^V \rightarrow (G/H)_+ \sma S^V$ and $\tau(\pi) \co (G/K)_+ \sma S^V \rightarrow (G/H)_+ 
\sma S^V$ (where $\pi$ is the projection $G/H \rightarrow G/K$, $H \subset K \subset G$) 
yield transfer homomorphisms in the cohomology theory, and these in turn give rise to 
transfers $M(G/H) \rightarrow M(G/K)$ and $M(G/H) \rightarrow M(G/G)$ respectively.  The 
construction of these transfer maps is intimately connected to Atiyah duality for the orbit 
spectra $\Sigma^\infty G/H_+$; this relationship can be seen via an explicit construction in 
terms of Pontryagin--Thom maps \cite[IX.3]{may-alaska} or from the perspective of formal 
categorical duality \cite[XVII.1]{may-alaska}.

There is a model theoretic interpretation of \fullref{t:char}, obtained by situating the 
theorem in the context of the study of diagram spectra.  The work of Lydakis \cite{lydakis} 
and Mandell, May, Schwede and Shipley \cite{mmss} permits the following modern 
reinterpretation of the $\Gamma$--space approach to infinite loop space theory.  Denote by 
$\sP$ the category of prespectra and $\sF`` \sT$ the category of $\Gamma$--spaces.  Let 
$\sW\! \sT$ be the category of continuous functors from $\sW$ to based spaces.  Our interest 
in $\sW\! \sT$ is its intermediate position between $\Gamma$--spaces and prespectra.  These 
various categories of spectra are linked by adjoint pairs $(\bP, \bU)$, where we regard the 
left adjoint $\bP$ as prolongation and the right adjoint $\bU$ is the restriction:
\[
\xymatrix{ \relax \sP \ar@<1ex>[r]^{\bP} & \ar@<1ex>[l]^{\bU} \relax\sW\! \sT 
\ar@<1ex>[r]^{\bU} & \ar@<1ex>[l]^{\bP} \relax \sF`` \sT}
\]
With suitable stable model structures, the first adjoint pair induces a Quillen equivalence 
and the second adjoint pair is a Quillen adjunction which induces an equivalence of the 
respective homotopy categories of connective objects (a ``connective'' Quillen equivalence 
\cite{mmss}).  From this standpoint, the work of Segal \cite{segal} amounts to the proof that 
for a ``very special'' $\Gamma$--space $E$, the $\sW$--space $\bP` E$ is fibrant and 
therefore the prespectrum $\bU \bP` E$ is fibrant and hence an $\Omega$--prespectrum.

To understand the situation equivariantly, it is natural to ask how much of this analysis can 
be generalized.  When $G$ is finite, one can obtain an identical version of this diagram.  
The homotopical analysis of equivariant $\Gamma$--spaces is known \cite{equi-segal,shimakawa} 
and we intend to discuss the model theoretic aspects of this elsewhere. When $G$ is an 
infinite compact Lie group, even though we cannot hope to have a ``connective'' equivalence 
of the equivariant analogues of $\sW$--spaces and $\Gamma$--spaces, a concrete understanding 
of the lefthand terms in the comparison diagram should indicate the nature of the 
generalization of $\Gamma$--spaces that will be needed.

We begin by carrying out an essentially formal reworking of the model category theory 
associated to $\sW$--spaces \cite{mmss} in the equivariant setting.  Let $G` \sW\! \sT$ 
denote the category of $\sWuG$--spaces with morphisms the $G$--equivariant natural 
transformations.  Fix a complete universe $U$ of real $G$--representations.

We have an adjoint pair of functors $(\bP, \bU)$ connecting the category of $G$--prespectra 
indexed on $U$ and the category of $\sWuG$--spaces, where the left adjoint $\bP$ is a 
prolongation and the right adjoint $\bU$ is a restriction:
\[
\xymatrix{ \relax G` \sP \ar@<1ex>[r]^{\bP} & \ar@<1ex>[l]^{\bU} \relax G` \sW\! \sT}
\]
This allows us to define the homotopy groups of a $\sWuG$--space $Z$ as the equivariant 
homotopy groups of the associated prespectrum $\bU Z$. That is, $\pi_*(Z) = \pi_* \bU Z$, 
where the homotopy groups run over fixed-point spaces corresponding to all closed subgroups 
of $G$ \cite[3.3.2]{mandell-may}.  In addition, $G` \sW\! \sT$ is symmetric monoidal with a 
smash product constructed by Kan extension \cite[2.3.1]{mandell-may}.  The comparison to 
$G$--prespectra factors through the category of orthogonal $G$--spectra (denoted $G`` \sI\! 
\sS$) via another adjoint pair $(\bP, \bU)$
\[
\xymatrix{ \relax G`` \sI\! \sS \ar@<1ex>[r]^{\bP} & \ar@<1ex>[l]^{\bU} \relax G` \sW\! \sT}
\]
where $\bP$ and $\bU$ are respectively strong and lax symmetric monoidal.  With this 
framework, we prove the following theorem recapitulating the nonequivariant theory of 
$\sW$--spaces \cite[17.1--17.6]{mmss} in the equivariant setting.  Note that a 
``$\!G$--topological'' model structure satisfies an appropriate analogue of Quillen's axiom 
SM7 reflecting compatibility with the $G$--enrichment \cite[3.1.4]{mandell-may}.

\begin{thm}
Fix a complete $G$--universe $U$. 
\begin{enumerate}
\item{There is a cofibrantly generated $G$--topological model
    structure on the category $G` \sW\! \sT$ in which the weak equivalences are the
    $\pi_*$--equivalences, the ``stable model structure''.}
\item{The fibrant objects in the stable model structure on $G` \sW\! \sT$
    are the $\sWuG$--spaces $Z$ such that $\bU Z$ is an
    $\Omega$--$G$--prespectrum.}
\item{The stable model structures satisfies the monoid and
    pushout-product axioms with respect to the smash product, and
    hence can be lifted to a model structure on categories of rings
    and modules.}
\item{The adjoint pair $(\bP, \bU)$ connecting $G` \sW\! \sT$ to
    $G`` \sI\! \sS$ is a Quillen equivalence.}
\item{There is a different model structure on $G` \sW\! \sT$ which is a
    cofibrantly generated $G$--topological model structure in which the
    weak equivalences are the $\pi_*$--equivalences, the
    ``absolute stable model structure''.  The identity functor is the left
    adjoint of a Quillen equivalence between the stable model
    structure and the absolute stable model structure.}
\item{The fibrant objects in the absolute stable model structure on $G
    \sW\! \sT$ are the $\sWuG$--spaces $Z$ such that the collection
    $\{Z(A \sma S^V)\}$ as $V$ varies over $U$ is an
    $\Omega$--$G$--prespectrum for any $A \in \sWuG$.}
\end{enumerate}
\end{thm}

By itself, this formal analysis has not bought us very much new information.  However, let us 
again recall the nonequivariant situation. There, one proves that the prolongation of a 
``very special'' $\Gamma$--space $X$ yields a $\sW$--space $\tilde{X}$ which is fibrant in 
the absolute stable model structure.  A $\sW$--space which is fibrant in the absolute stable 
model structure is clearly fibrant in the stable model structure as well, and then the 
restriction $\bU \tilde{X}$ is a fibrant prespectrum.  It is important to be clear about the 
roles of the two stable model structures on $\sW$--spaces. Although at the end of the day we 
are interested in fibrant objects in the stable model structure on $\sW$--spaces, the 
absolute stable model structure is essential in order to compare $\sW$--spaces to 
$\Gamma$--spaces.  In addition, it turns out that the fibrant objects in the absolute stable 
model structure admit a concise intrinsic description in the nonequivariant setting---they 
are precisely the linear functors.

Thus, we are led to the question of determining a similar intrinsic description of the 
$\sWuG$--spaces $Z$ which are fibrant in the absolute stable model structure.  Explicitly, we 
want conditions which are necessary and sufficient for the prespectra $\{Z(A \sma S^V)\}$ to 
be $\Omega$--$G$--prespectra.  Such conditions are provided by \fullref{t:char}.
\vspace{-6pt}

\begin{corollary}
A $\sWuG$--space is fibrant in the absolute stable model structure on $G` \sW\! \sT$ if and 
only if it is equivariantly linear.
\end{corollary}
\vspace{-6pt}

The concrete characterization of equivariantly linear functors provides an indication of the 
information that must be captured by a space-level recognition principle for equivariant 
infinite loop spaces.  For instance, it suggests that the appropriate equivariant analogue of 
$\Gamma$--spaces involve an enlarged domain category which contains all the orbit spaces 
$G/H$.  In future work, we intend to exploit this perspective.
\vspace{-6pt}

The paper is organized as follows.  In section 2, we briefly state the definitions and model 
theoretic results we will refer to in the course of proving the main theorem.  We relegate 
proofs to the appendix.  In section 3, we prove \fullref{t:char}.  In the first section of 
the appendix, we carry out the model theoretic analysis of the category of $\sWuG$--spaces.  
This is very similar to the analysis of \cite{mmss}, and our primary purpose is to record 
results and the proofs of supporting lemmas which do not follow immediately from the 
nonequivariant results.  In the second part of the appendix, we analyze 
$\Gamma$--$S^1$--spaces.
\vspace{-6pt}

\paragraph*{Acknowledgements} I wish to express my gratitude to Peter May for proposing 
this area of research, for many useful discussions, and for careful reading of previous 
drafts.  I would also like to thank Mike Mandell for his invaluable assistance and 
suggestions.  Brooke Shipley and Chris Douglas provided valuable comments based on an earlier 
draft, and I had helpful conversations with Mark Behrens, Ben Lee, and Ben Wieland. In 
particular, I would like to thank Ben Wieland for catching a mistake in a previous version of 
the paper.  The exposition of this paper was significantly improved by comments from an 
anonymous referee.

\section{Basic definitions and a rapid overview of model theoretic
  results}
\vspace{-6pt}

In this section, we will present the basic definitions and summarize
the model theoretic results.  The proofs of these results appears in
the first section of the appendix.
\vspace{-6pt}

\subsection[Categories of S_G-spaces]{Categories of $\sWuG$--spaces}
\vspace{-6pt}

The categories we will be working with are enriched over based $G$--spaces.  Thus, our 
discussion could be cast entirely in terms of enriched category theory.  However, we follow 
the convention of Mandell and May \cite{mandell-may} and instead consider ordinary functors 
with additional conditions in order to emphasize the analogies to the nonequivariant case and 
to minimize overhead.  As a consequence, we will work with pairs $\sC_G$ and $G`\sC = 
(\sC_G)^G$, where $\sC_G$ is a category of $G$--objects and nonequivariant maps between them, 
and $G`\sC$ is obtained by restricting to the $G$--maps.  The hom spaces of $\sC_G$ are given 
a $G$--action via conjugation, and regarded as based via the addition of a $G$--fixed 
basepoint as necessary. Therefore, we can obtain the space of $G$--maps by taking fixed 
points. For instance, let $\sT_G$ be the category of based $G$--spaces with morphisms all 
maps of nonequivariant spaces.  Then $\hom_{\sT_G}(X,Y)$ has a $G$--action given by 
conjugation.  The category $G`\sT$ is obtained by passage to $G$--fixed points on the hom 
spaces, and is the category of $G$--spaces and $G$--maps.  All of the model structures we 
consider are compatible with this enrichment; we will refer to such model categories as 
$G$--topological.  The precise definition of this compatibility is discussed in the appendix.
\vspace{-6pt}

With this in mind, we can now define the categories of $\sWuG$--spaces we will be working 
with.  Recall that $\sWuG$ denotes the category of based $G$--spaces homeomorphic to finite 
$G$--CW--complexes, with morphisms all continuous (but not necessarily equivariant) maps.  
This is a full subcategory of $\sT_G$.  Although $\sWuG$ is not small, it is skeletally small 
and throughout we will tacitly assume that we have chosen a small skeleton of $\sWuG$ and are 
working relative to that skeleton.
\vspace{-6pt}

\begin{definition}
The category $\sWuG \sT$ has as objects the based continuous $G$--functors from $\sWuG$ to 
$\sT_G$.  The morphisms are the natural transformations between functors.  We can topologize 
the morphisms as a subspace of the product of the function spaces $\Map_0(X(A),Y(A))$ over $A 
\in \sWuG$, which has a $G$--action by conjugation.
\vspace{-6pt}

The category $G`\sW\!\sT$ of $\sWuG$--spaces is obtained by passage to $G$--fixed points from 
$\sWuG \sT$.  That is, the objects are again the continuous $G$--functors from $\sWuG$ to 
$\sT_G$ and the morphisms are the natural $G$--maps.
\end{definition}
\vspace{-6pt}

Just as in the case of orthogonal $G$--spectra \cite[2.3.1]{mandell-may}, we have the 
following theorem.  

\begin{thm}
The categories $\sWuG\sT$ and $G`\sW\!\sT$ have smash product and function spectrum functors 
which make them closed symmetric monoidal categories.  The unit is the identity functor.
\end{thm}

The smash product is constructed in the usual way as an
internalization of the obvious external smash product via left Kan
extension.  Some care has to be taken to verify that the Kan extension
exists \cite[2.6.7]{mandell-may}.

\subsection{The stable model structure}

Throughout, fix a complete universe $U$.  The first model structure we consider is the 
relative level model structure on $G`\sW\! \sT$, where by relative we mean that the 
fibrations and weak equivalences are detected only on the spheres $\{S^V\}$ for $V \in U$.

\begin{definition}
The relative level model structure on this category is defined as
follows.  A map $Y \rightarrow Z$ is
\begin{enumerate}
\item{a fibration if each $Y(S^V) \rightarrow Z(S^V)$ is an
equivariant Serre fibration,}
\item{a weak equivalence if each $Y(S^V) \rightarrow Z(S^V)$ is an
equivariant weak equivalence,}
\item{a cofibration if it has the left-lifting property with respect
to the acyclic fibrations.}
\end{enumerate}
\end{definition}

\begin{proposition}
The relative level model structure on $G`\sW\!\sT$ is a cofibrantly generated 
$G$--topological model structure.
\end{proposition}

There is an associated stable model structure.  We define $\pi_* Z$ for a $\sWuG$--space $Z$ 
by passing to the $G$--prespectrum $\bU Z$ and specifying $\pi_* Z = \pi_* \bU Z$.  Recall 
that for a subgroup $H$ of $G$ and an integer $q$, for $q \geq 0$ we define
\[\pi^H_q(\bU Z) = \colim_V \pi_q^H (\Omega^V Z(V))\]
and for $q > 0$ we define
\[\pi^H_{-q}(\bU Z) = \colim_{V \supseteq \R^q} \pi_0^H (\Omega^{V -
  \R^q} Z(V)).\]
These equivariant homotopy groups capture stable equivalences.

\begin{definition}
Let $Z_1$ and $Z_2$ be $\sWuG$--spaces.  A map $f\co X_1 \rightarrow X_2$ is a 
$\pi_*$--isomorphism if the induced maps $f_* \co \pi^H_q(X_1) \rightarrow \pi^H_q(X_2)$ are 
isomorphisms for all closed subgroups $H \subset G$.
\end{definition}
\eject
\begin{definition}
In the stable model structure, a map is 
\begin{enumerate}
\item{a cofibration if it is a cofibration in the relative level model structure,}
\item{a weak equivalence if it is a $\pi_*$--equivalence,}
\item{a fibration if it has the right-lifting property with respect to the acyclic cofibrations (maps which are both level cofibrations and $\pi_*$--equivalences).}
\end{enumerate}
\end{definition}

\begin{proposition}
The stable model structure is a $G$--topological model structure on the category $G` \sW\! 
\sT$.
\end{proposition}

There is a pair of adjoint functors ($\bP$, $\bU$) connecting $G` \sW\! \sT$ and the category 
$G``\sI\! \sS$ of orthogonal $G$--spectra \cite{mandell-may,mmss}.  $\bU$ is the forgetful 
functor from $\sWuG$--spaces to orthogonal $G$--spectra, and $\bP$ is the prolongation 
constructed as a left Kan extension along the inclusion of domain categories 
\cite[23.1]{mmss}.  We have the expected comparison result.

\begin{thm}
The pair $(\bP,\bU)$ specifies a Quillen equivalence between the stable model category 
structure on $G` \sW\! \sT$ and the stable model category structure on $G`` \sI\! \sS$.
\end{thm}

\subsection{The absolute stable model structure}

The level model structure used to construct the stable model structure in the previous 
section depends on evaluation at the spheres.  This makes it inconvenient to compare to 
diagram categories where the domain does not include an embedding of the spheres, for 
instance $\Gamma$--$G$--spaces.  As in the nonequivariant case, we rectify this by 
considering an ``absolute'' model structure.

\begin{definition}
The absolute level model structure on $G` \sW\! \sT$ is defined as follows.  A map $Y 
\rightarrow Z$ is
\begin{enumerate}
\item{a fibration if each $Y(A) \rightarrow Z(A)$ for $A \in \sWuG$ is
  an equivariant Serre fibration,}
\item{a weak equivalence if each $Y(A) \rightarrow Z(A)$ for $A \in
  \sWuG$ is an equivariant weak equivalence,}
\item{a cofibration if it has the left-lifting property with respect
  to the acyclic fibrations.}
\end{enumerate}
\end{definition}

There is an associated absolute stable model structure.
\eject
\begin{definition}
In the absolute stable model structure, a map is
\begin{enumerate}
\item{a cofibration if it is a cofibration in the absolute level model structure,}
\item{a weak equivalence if it is a $\pi_*$--equivalence,}
\item{a fibration if it has the right-lifting property with respect to the acyclic cofibrations.}
\end{enumerate}
\end{definition}

\begin{proposition}
The stable model structure is a $G$--topological model structure on the category $G` \sW\! 
\sT$.
\end{proposition}

In the course of proving the previous proposition, we obtain an
identification of the fibrant objects in the absolute stable model
structure.

\begin{proposition}
A $\sWuG$--space $Z$ is fibrant in the absolute stable model structure if and only if for all 
$A \in \sWuG$ and $W \in U$ the structure map 
\[Z(A) \rightarrow \Omega^W Z(S^W \sma A)\] 
is a weak equivalence.
\end{proposition}

There is a Quillen equivalence between the relative stable model
structure and the absolute stable model structure.

\begin{thm}
The identify functor is the left adjoint of a Quillen equivalence between the category of 
$\sWuG$--spaces with the relative stable model structure and the category of $\sWuG$--spaces 
with the absolute stable model structure. 
\end{thm}

\subsection{Ring and module spectra}

Just as in the nonequivariant setting, we can lift the stable model structure to categories 
of ring and module $\sWuG$--spaces.  While we do not require this in the paper, it is a 
useful feature of this perspective on the equivariant stable category.

\begin{thm}
Let $R$ be a ring $\sWuG$--space.
\begin{enumerate}
\item{The category of $R$--module $\sWuG$--spaces is a cofibrantly generated proper
$G$--topological model category, with weak equivalences and fibrations created in the stable 
model structure on the category of $\sWuG$--spaces.}
\item{If $R$ is commutative, the category of $R$--algebra
$\sWuG$--spaces is a cofibrantly generated right proper $G$--topological model category with 
weak equivalences and fibrations created in the stable model structure on the category of 
$\sWuG$--spaces.}
\end{enumerate}
\end{thm}

\begin{remark}
The obvious variant of this theorem starting from the absolute stable
model structure also holds.
\end{remark}

\section[Continuous G-functors and Omega-G-prespectra]{Continuous $G$--functors and $\Omega$--$G$--prespectra}

In this section we will provide concrete conditions which describe the fibrant objects in the 
absolute stable model structure on $\sWuG$--spaces.  That is, we specify conditions on a 
$\sWuG$--space $Z$ which guarantee that for any $A \in \sWuG$, the prespectrum obtained as 
the collection of spaces $\{Z(S^V \sma A)\}$ is a genuine $\Omega$--$G$--spectrum.  These 
conditions amount to enforcing a suitable interaction with equivariant Spanier--Whitehead 
duality (or more precisely equivariant Atiyah duality) for orbit spectra.  This connection to 
duality highlights the difficulty of generalizing recognition principles from the case of $G$ 
finite to the case of $G$ a compact Lie group, for only when $G$ is finite are the orbit 
spectra self dual.

\subsection{Linearity in the nonequivariant setting}

First, we recall the nonequivariant situation.  Let $Z$ be a
continuous functor from $\sW$ to spaces.  There is a structure map
$\sigma \co Z(A) \sma B \rightarrow Z(A \sma B)$ \cite[4.9]{mmss} which
arises as a consequence of continuity.  The map $\sigma$ is the
adjoint of the composite
\[
\xymatrix{ B \ar[r]^-\alpha & \relax\sT(A, A \sma B) \ar[r]^-= & \relax \sW(A, A \sma B) 
\ar[r]^-Z & \relax\sT(X(A), X(A \sma B)), }\] where $\alpha(b)(a) = a \sma b$.  Setting $B = 
S^n$, this gives us the structure maps of a prespectrum when we consider the collection 
$\{Z(S^n)\}$.  We will denote this prespectrum by $Z[S^0]$, and write $Z[A]$ for the 
prespectrum $\{Z(S^n \sma A)\}$.  Observe that $Z[S^0] = \bU Z$.  Also, note that the 
structure maps $Z(A) \sma I_+ \rightarrow Z(A \sma I_+)$ imply that any $\sW$--space 
preserves homotopies and hence weak equivalences on $\sW$ \cite[17.4]{mmss}.

It is well known that there is a simple condition on $Z$ which guarantees that $Z[A]$ is an 
$\Omega$--prespectrum.  For this to be true, it must be the case that $Z$ takes homotopy 
pushout squares to homotopy pullback squares.  This is sometimes stated as $Z$ is linear. 
Given such a $Z$, for any $A \in \sW$ we have the homotopy pushout 
\[
\begin{CD}
A @>>> * \\
@VVV @VVV \\
* @>>> S^1 \sma A \\
\end{CD}
\]
which constructs the suspension, and when we apply $Z$ to this diagram
there is an induced weak equivalence between $Z(A)$ and $\Omega
Z(\Sigma A)$.

We now wish to generalize this to describe similar conditions in the case of a $\sWuG$--space 
$Z$ which will guarantee that the collection $\{Z(S^V \sma A)\}$ forms an 
$\Omega$--$G$--spectrum.  Denote the prespectrum $\{Z(S^V \sma A)\}$ by $Z_U[A]$.  Again, 
note that $Z_U[S^0] = \bU Z$.  The structure maps arise via an adjunction analogous to the 
nonequivariant case: 
\[
\xymatrix{ B \ar[r]^-\alpha & \relax\sT_G(A, A \sma B) \ar[r]^-= & \relax \sWuG(A, A \sma B) 
\ar[r]^-Z & \relax\sT_G(Z(A), Z(A \sma B)) }\] Once again, the existence of these structure 
maps implies that $\sWuG$--spaces preserve weak equivalences in $\sWuG$.  Next, observe that 
taking homotopy pushout squares to homotopy pullback squares (in the category of $G$--spaces) 
is insufficient to handle desuspension by arbitrary representations.  We cannot construct 
$\Omega^V Z(S^V)$ for arbitrary representations $V$ in the fashion above.  This condition is 
however enough to construct a naive $\Omega$--$G$ spectrum (indexed on a trivial universe).

Following \cite[17.9]{mmss}, we obtain the following characterization of continuous functors 
which generate naive $\Omega$--$G$--prespectra.

\begin{proposition}\label{P:naive}
Let $Z$ be a $\sWuG$--space.  The following are equivalent:
\begin{enumerate}
\item{$Z$ takes $G$--homotopy pushout squares to $G$--homotopy pullback squares.}
\item{For any $A \in \sWuG$, the prespectrum $Z[A]$ is a naive $\Omega$--$G$--prespectrum.}
\item{For any $A \in \sWuG$, the adjoint structure map $Z(A)
  \rightarrow \Omega Z(\Sigma A)$ is an equivalence.}
\end{enumerate}
\end{proposition}

We will call such $\sWuG$--spaces ``naively equivariantly linear''.  In order to handle 
suspensions at arbitrary representations, we need to specify more data about the functor $Z$.

\subsection{Compatibility with equivariant Spanier--Whitehead duality}

Fix a complete universe $U$.  Let $V$ be a $G$--representation in $U$ and let $G/H$ be an 
orbit $G$--space which is embedded in $V$.  Denote by $T\nu$ the Thom space of the normal 
bundle $\nu$ of the embedding, and note that this is $G$--homeomorphic to the 
compactification of a tubular neighborhood of $G/H$ because $G/H$ is compact.  More 
concretely, for sufficiently small $\epsilon$ we can describe the tubular neighborhood of 
$G/H$ in $V$ as $(G/H)_\epsilon$, the $\epsilon$--neighborhood of $G/H$.  Then $T\nu$ is 
$G$--homeomorphic to $(G/H)_\epsilon^c$, the one-point compactification of $(G/H)_\epsilon$.  
We know that $G/H_+$ and $(G/H)_\epsilon^c$ are equivariantly $V$--dual.

This duality can be exhibited by the following map.  There is a map $$G/H \rightarrow 
\Map(D(\epsilon), (G/H)_\epsilon)$$ taking each $m \in G/H$ to the map which takes an element 
$x$ of the $\epsilon$--ball $D(\epsilon)$ about the origin to $m + x$.  This induces a based 
map 
\[(G/H)_+ \rightarrow \Map_0((G/H)_\epsilon^c, D(\epsilon)^c) \cong
\Map_0((G/H)_\epsilon^c, S^V)\] by taking an element of
$(G/H)_\epsilon^c$ to the basepoint if it is not within $\epsilon$ of
$m$ and to $x - m$ otherwise.  The adjoint of this is the duality map 
\[G/H_\epsilon^c \sma (G/H)_+ \rightarrow S^V.\]
This is the classical Atiyah duality map.  Note that we don't actually use the fact that 
$G/H$ is a submanifold of $V$.  This map makes sense whenever we have a compact $G$--subset 
of $V$, and is the duality map when $\epsilon$ is sufficiently small.

Now assume that we have a $\sWuG$--space $Z$.  Given the map 
\[(G/H)_+ \rightarrow \Map_0((G/H)_\epsilon^c, S^V)\]
by functoriality we obtain a map
\[(G/H)_+ \rightarrow \Map_0(Z((G/H)_\epsilon^c), Z(S^V))\]
and by adjunction we have a map
\[\xi \co Z((G/H)_\epsilon^c) \rightarrow \Map_0((G/H)_+, Z(S^V)).\]
Now, if we have a space $X \in \sWuG$ and smash the duality map on both sides by $X$, by the 
same process we obtain a map
\[\xi \co Z((G/H)_\epsilon^c \sma X) \rightarrow \Map_0((G/H)_+, Z(S^V
\sma X)).\]

\begin{remark}\label{r:xi}
One subtlety of the duality theory developed by Lewis, May and Steinberger 
\cite[3.1--3.8]{lms} is that for given $V$--duals $X$ and $Y$ there are many possible choices 
of space-level maps exhibiting the $V$--duality.  In the specific case of $G/H$ (and more 
generally for embedded submanifolds), there is another very explicit description of the 
duality between $G/H_+$ and $T\nu$, the Thom space of the normal bundle of the embedding.

Specifically, we can construct a map
\[T\nu \sma G/H_+ \rightarrow G/H_+ \sma S^V \rightarrow S^V\]
where the first map is a Pontryagin--Thom map associated with a tubular neighborhood of the 
composite \[G/H \rightarrow G/H \times G/H \rightarrow \nu \times G/H\] and the second is the 
collapse map of $G/H_+$ onto $S^0$ \cite[3.5.1]{lms}.  By functoriality and manipulation of 
adjoints we can obtain a map 
\[\xi_2 \co Z(T\nu) \rightarrow \Map_0(G/H_+, Z(S^V))\]
analogously to the construction of $\xi$.

Under the homeomorphism given by the tubular neighborhood theorem
$T\nu \rightarrow G/H_\epsilon^c$, the maps $\xi$ and $\xi_2$
coincide \cite[3.5.1]{lms}.  As a consequence, our conditions below
could be phrased in terms of the map $\xi_2$ rather than $\xi$, and in
general, we could phrase them abstractly in terms of any system of
duality maps which are suitably functorial.
\end{remark}

We are now ready to formulate the first version of the additional conditions required for a 
$\sWuG$--space to represent a genuine $\Omega$--$G$--prespectrum.  We refer to the condition 
below as ``hypothesis (A)''.

\begin{hyp}\label{d:xi}
A $\sWuG$--space $Z$ satisfies \textit{hypothesis (A)} for the universe $U$ if the following 
two conditions hold.
\begin{enumerate}
\item{$Z$ takes $G$--homotopy pushout squares to $G$--homotopy pullback squares.}
\item{For all $X \in \sWuG$ and any smooth embedding $G/H \monoto V$,
    $V \in U$,
    there exists $\epsilon > 0$ such that the map 
\[\xi \co Z((G/H)_\epsilon^c) \sma X) \rightarrow \Map_0(G/H_+, Z(S^V \sma X))\] 
is a $G$--equivalence.}
\end{enumerate}
\end{hyp}

We will need a lemma extending this condition slightly.

\begin{lemma}\label{l:easy}
If $Z$ satisfies hypothesis (A) and $G/H$ embeds in $U$, then for any
smooth embedding $G/H \times D^n \monoto V$, there exists
$\epsilon > 0$
 such that the map
\[\xi \co Z((G/H \times D^n)_\epsilon^c) \sma X) \rightarrow \Map_0(G/H_+, Z(S^V
\sma X))\] is a $G$--equivalence.
\end{lemma}

First, we show that this condition is sufficient.  The argument below
is an adaptation of the argument due to Segal \cite{equi-segal} (and
corrected by Shimakawa \cite{shimakawa}) for the case when $G$ is a
finite group.  Recall that we wrote $Z_U[A]$ to refer to the
prespectrum $\{Z(S^V \sma A)\}$.

\begin{thm} \label{T:segal}
Let $Z$ be a $\sWuG$--space which satisfies hypothesis (A) for the universe $U$.  Then 
$Z_U[A]$ is an $\Omega$--$G$--prespectrum.
\end{thm}

\begin{proof}
Fix an arbitrary representation $V$.  By naive linearity, without loss of generality we can 
assume that $V$ contains a trivial representation $\R$.  This assumption allows us to provide 
$G$--fixed basepoints to subspaces of $V$.  Let $D(1)$ denote the unit disk in $V$, and 
$S(1)$ the unit sphere which is the boundary, and in general let $D(r)$ and $S(r)$ be the 
disk and its boundary of radius $r$.  There is a commutative diagram 
\[
\begin{CD}
Z(D(1+\epsilon)/S(1+\epsilon)) @>>> \Map_0(D(1)_+, Z(S^V)) \\
@VVV @VVV \\
Z(D(1+\epsilon) / (D(1-\epsilon) \cup S(1+\epsilon))) @>>> \Map_0(S(1)_+, Z(S^V)) \\
\end{CD}
\]
where the horizontal maps are defined analogously to the map $\xi$.

The top horizontal map is clearly an equivalence, and temporarily
assume the bottom horizontal map is an equivalence.  The right
vertical map is a fibration induced from the inclusion $S(1)_+
\monoto D(1)_+$, and the fiber is $\Omega^V Z(S^V)$.  There is a
cofibration sequence 
\[(D(1-\epsilon) \cup S(1+\epsilon)) / S(1+\epsilon)
\rightarrow D(1+\epsilon) / S(1+\epsilon) \rightarrow D(1+\epsilon) / (D(1 - \epsilon) \cup 
S(1+\epsilon)).\] Since $(D(1-\epsilon) \cup S(1+\epsilon)) / S(1+\epsilon) \htp S^0$ and $Z$ 
takes cofibration sequences to fibration sequences by hypothesis, we know that the homotopy 
fiber of the left vertical map is $Z(S^0)$.  Therefore, we can conclude that there is an 
equivalence $Z(S^0) \simeq \Omega^V Z(S^V)$.  The induced map of fibers from $Z(S^0) 
\rightarrow \Omega^V Z(S^V)$ is indeed the adjoint of the structure map, as it is obtained 
from the embedding of $0$ in $V$. Since $Z(A \sma S^W \sma -)$ is also a functor satisfying 
our hypotheses, we obtain the desired equivalence $Z(S^W \sma A) \simeq \Omega^V Z(S^{V\oplus 
W} \sma A)$ for every $A \in \sWuG$.

Therefore, to complete the proof of the theorem it will suffice to
verify that the bottom map is indeed a weak equivalence.  As an aside,
note that the bottom map can be described as the map $\xi_2 \co Z(T\nu)
\rightarrow \Map_0(S(1)_+, Z(S^V))$, where $T\nu$ here is Thom space of
the normal bundle of the embedding of $S(1)$ in $V$.

Since $S(1)$ is a finite $G$--CW--complex, we proceed by induction.  We can decompose $S(1)$ 
as a regular $G$--CW--complex \cite{illman}.  That is, we can regard it as comprised of cells 
$G/H \times D^n$ (for varying $H$) where the attaching maps are homeomorphisms and the images 
of the boundary $G/H \times S^{n-1}$ are equal to unions of cells of lower dimension.  
Moreover, by subdividing if necessary, we can assume that the closed cells $G/H \times D^n$ 
are subcomplexes. We will fix a choice of homeomorphic embedding of $S(1)$ in $V$. Recall 
that $D^n$ has trivial $G$--action.

Let $X$ be a subcomplex of $S(1)$, a union of some of the cells of $S(1)$.  We
have maps 
\[\xi_X \co Z(X^c_\epsilon) \rightarrow \Map_0(X_+, Z(S^V))\]
which we can regard as induced by restriction of the map \[S(1)
\rightarrow \Map_0(D(\epsilon), S(1)_\epsilon)\]
\newpage
used to construct $\xi$.
\medskip

We will induct downward over the number of cells in $X$.  The base
cases therefore involve $X$ consisting of a single cell $G/H \times
D^n$.  In this situation, the map $\xi_X$ is an equivalence by \fullref{l:easy}.
\medskip

Let the number of cells in $X$ be $m$, and assume that $\xi$ is an
equivalence for subcomplexes with $m-1$ cells or fewer.  Let $G/H
\times D^n$ be a cell of highest dimension in $X$, and let $Y$ be the
union of the remaining cells, so that $X = Y \cup (G/H \times D^n)$.
\medskip

There is a commutative diagram
\[
\xymatrix@C=13pt{ ((Y {\cup} (G/H {\times} D^n))_\epsilon {-} (G/H {\times} D^n)_\epsilon)^c
  \ar[r] \ar[d] & (Y {\cup} (G/H {\times} D^n))_\epsilon^c \ar[r] \ar[d] & (G/H
{\times} D^n)_\epsilon^c \ar[d] \\
(Y_\epsilon - (Y \cap (G/H {\times} D^n))_\epsilon)^c \ar[r] &
Y_\epsilon^c \ar[r] & (Y {\cap} (G/H {\times} D^n))_\epsilon^c \\
}
\]
where each row is a cofibration.  The map 
\[((Y \cup (G/H \times D^n))_\epsilon - (G/H \times D^n)_\epsilon)^c
\rightarrow (Y_\epsilon - (Y \cap (G/H \times D^n))_\epsilon)^c\] is a weak equivalence, as 
follows.  Recall that for $\epsilon$ sufficiently small we can naturally replace the diagram 
above with the corresponding diagram of cones \cite[2.4.13]{lms}:
\[
\xymatrix@C=16pt{
C(V{-}(G/H {\times} D^n), V{-}X) \ar[r] \ar[d] & C(V,V{-}X) \ar[r] \ar[d] & C(V,V{-}(G/H {\times} D^n)) \ar[d] \\
C(V{-}((G/H {\times} D^n) {\cap} Y), V{-}Y) \ar[r] & C(V,V{-}Y) \ar[r] & C(V,V{-}(Y {\cap} (G/H {\times} 
D^n))) }
\]

Now consider the set $U = (G/H \times D^n) - ((G/H \times D^n) \cap
Y)$.  The closure of $U$ in $V- ((G/H \times D^n) \cap Y)$ is
contained in the interior of $V-Y$, and so excision \cite[2.4.3]{lms}
implies that the leftmost map is a weak equivalence.
\medskip

Therefore, upon application of $Z$ we obtain a homotopy pullback square:
\[
\xymatrix{
Z((Y \cup (G/H \times D^n))_\epsilon^c) \ar[r] \ar[d] & Z((G/H \times D^n)_\epsilon^c) \ar[d] \\
Z(Y_\epsilon^c) \ar[r] & Z((Y \cap (G/H \times D^n))_\epsilon^c) }
\]

In addition, applying the mapping space functor $\Map_0(-,Z(S^V))$ and
recalling that it also takes cofibrations to fibrations, we have a
homotopy pullback square 
\[
\xymatrix{
\Map_0((Y \cup (G/H \times D^n))_+, Z(S^V)) \ar[r]\ar[d] & \Map_0((G/H \times
D^n)_+, Z(S^V)) \ar[d] \\
\Map_0(Y_+, Z(S^V)) \ar[r] & \Map_0((Y \cap (G/H \times D^n))_+, Z(S^V)) \\
}
\]
and one checks that the cube induced by the maps $\xi_X$ is commutative.
Since $S(1)$ is regular, $Y \cap (G/H \times D^n)$ is a subcomplex
consisting of strictly fewer cells (of lower dimension) and so the
inductive hypothesis implies that we have equivalences at the three
nonterminal corners of the cube, and therefore there is an equivalence
at the terminal corner.
\end{proof}

We now wish to show that in fact hypothesis (A) is necessary.  In
order to do so, we must first recall the following equivariant version
of a theorem of Lydakis \cite{lydakis}.

\begin{proposition}\label{p:lyd}
Let $Z$ be a $\sWuG$--space and $A \in \sWuG$.  Then the maps $Z_U(X) \sma A \rightarrow 
Z_U(X \sma A)$ induce a $\pi_*$--equivalence $Z_U[S^0] \sma A \htp Z_U[A]$. 
\end{proposition}

\begin{proof}
The proof follows the nonequivariant proof given in \cite[17.6]{mmss}. In order to perform 
the induction, we substitute the equivariant theorems \cite[3.3.5]{mandell-may} for 
\cite[7.4]{mmss}.  Note that we depend on the fact that all $A \in \sWuG$ are stably 
dualizable when $U$ is a complete universe (as we are assuming here).  See 
\fullref{s:incomplete} for discussion of the situation when $U$ is not complete. 
\end{proof}

This result allows us to pass between $\pi_*$--equivalences of fibrant $\sWuG$--spaces and 
weak equivalences of the ``zero spaces'' of the $\sWuG$--spaces.  We will employ this 
observation to deduce information about our space-level maps $\xi$ from stable dualities.

\begin{notn}
For a prespectrum $Z$, we will sometimes write $\Omega^\infty Z$ in place of $Z(S^0)$.
\end{notn}

\begin{corollary} \label{C:lyd}
Let $Y$ and $Z$ be $\sWuG$--spaces which are fibrant in the absolute stable model structure.  
Denote the fibrant replacement of a prespectrum $D$ by $fD$.
\begin{enumerate}
\item{The spaces $\Omega^\infty f(Z_U[S^0] \sma A)$ and $Z(A)$ are weakly equivalent.}
\item{Given a map of prespectra $Y_U[S^0] \rightarrow Z_U[S^0]$
  induced from a natural transformation $Y \rightarrow Z$ and a
  map $A \rightarrow B$ for $A, B \in \sWuG$, if the induced map 
\[Y_U[S^0] \sma A \rightarrow Z_U[S^0] \sma B\] 
is a $\pi_*$--equivalence then the induced map
 \[Y(A) \rightarrow Z(B)\] is a weak equivalence.}
\item{Given a map of prespectra $Y_U[S^0] \rightarrow Z_U[S^0]$
  arising from a natural transformation $Y \rightarrow Z$ and a
  map $A \sma B \rightarrow C$ for $A, B, C \in \sWuG$, if the induced map 
\[Y_U[S^0] \sma A \rightarrow F(B, Z_U[S^0] \sma C)\]
is a $\pi_*$--equivalence then the induced map
\[Y(A) \rightarrow \Map_0(B, Z(C))\]
is a weak equivalence.}
\end{enumerate}
\end{corollary}

\begin{proof}
The first part is an immediate consequence of the preceding proposition.  Take any $A \in 
\sWuG$.  Since there is a $\pi_*$--equivalence $Z_U[S^0] \sma A \rightarrow Z[A]$, there is a 
$\pi_*$--equivalence $f(Z_U[S^0] \sma A) \rightarrow f(Z[A])$.  As these are 
$\Omega$--$G$--prespectra, a $\pi_*$--equivalence is the same as a level equivalence 
\cite[3.3.4]{mandell-may}, and so there is a weak equivalence $\Omega^\infty f(Z_U[S^0] \sma 
A) \rightarrow \Omega^\infty f(Z[A])$.  But since there is also a level equivalence $Z[A] 
\rightarrow f(Z[A])$, the result follows.

For the second claim, the given maps induce a commutative diagram 
\[
\begin{CD}
Y_U[S^0] \sma A @>>> Z_U[S^0] \sma B \\
@VVV @VVV \\
Y_U[A] @>>> Z_U[B] \\
\end{CD}
\] 
in which the vertical maps are $\pi_*$--equivalences.  Therefore, if the top horizontal map 
is a $\pi_*$--equivalence, the bottom map must also be a $\pi_*$--equivalence.  Since 
$Y_U[A]$ and $Z_U[B]$ are $\Omega$--$G$--prespectra by hypothesis, a $\pi_*$--equivalence is 
a level equivalence and therefore we have a weak equivalence $Y(A) = \Omega^\infty Y_U[A] 
\rightarrow \Omega^\infty Z_U[A] = Z(A)$.  Using the naturality of the structure maps, we can 
see that this induced weak equivalence coincides with the map $Y(A) \rightarrow Z(B)$ induced 
from the natural transformation $Y \rightarrow Z$ and the map $A \rightarrow B$.

Finally, the last part follows from an argument similar to the second
part.  There is a commutative diagram
\[
\begin{CD}
Y_U[S^0] \sma A @>>> F(B, Z_U[S^0] \sma C) \\
@VVV @VVV \\
Y_U[A] @>>> F(B, Z_U[C])\\
\end{CD}
\]
which arises as the adjoint of the commutative diagram
\[
\begin{CD}
Y_U[S^0] \sma A \sma B @>>> Z_U[S^0] \sma C \\
@VVV @VVV \\
Y_U[A] \sma B @>>> Z_U[C] \\
\end{CD}
\]
where the bottom vertical map is the composite of the structure map
and the given map $A \sma B \rightarrow C$.  Now we argue as above,
using the fact that for any space $B$ and prespectrum $Z$,
$\Omega^\infty F(B, Z) = \Map_0(B, \Omega^\infty Z)$. 
\end{proof}

With this in hand, we can complete the proof that our condition on $\sWuG$--spaces is 
necessary and sufficient for the prespectra $Z_U[A]$ to be $\Omega$--$G$--prespectra.

\begin{thm}\label{t:equilin}
For a $\sWuG$--space $Z$, the following are equivalent.
\begin{enumerate}
\item{$Z$ satisfies hypothesis (A).}
\item{For any finite $G$--CW--complex $A$, $Z_U[A]$ forms an
    $\Omega$--$G$--prespectrum.}
\item{For any finite $G$--CW--complex $A$, for any $W$ in
    the universe $U$, the adjoint of the
    structure map $Z(A) \sma S^W \rightarrow Z(S^W \sma A)$ is a weak
    equivalence $$Z(A) \simeq \Omega^W Z(\Sigma^W A).$$}
\end{enumerate}
\end{thm}

\begin{proof}
We have already shown that if $Z$ satisfies hypothesis $(A)$, then $Z_U[A]$ forms an 
$\Omega$--$G$--prespectrum for all $A \in \sWuG$.  By definition this is equivalent to the 
third condition.  Now assume that $Z[A]$ is an $\Omega$--$G$--prespectrum for every $A \in 
\sWuG$.  The ``naive'' version of this result, \fullref{P:naive}, implies that $Z$ takes 
$G$--homotopy pushout squares to $G$--homotopy pullback squares.

Thus, we need to show that the map $\xi$ is an equivalence for all $G/H$ which embed in $V$ 
for any $V$ in the universe $U$.  Let $E$ denote $Z_U[S^0]$, and recall this is an 
$\Omega$--$G$ prespectrum by hypothesis.  By the third part of \fullref{C:lyd}, 
\[\xi \co Z(T\nu \sma X) \rightarrow \Map_0(G/H_+, Z(S^V \sma X))\]
will be a weak equivalence if the map
\[\tilde{\xi} \co E \sma X \sma (G/H)_\epsilon^c \rightarrow F(G/H_+, E \sma S^V \sma
X)\] is a $\pi_*$--equivalence.  Here $\tilde{\xi}$ is obtained by adjunction from the map
\[
\xymatrix{ E \sma X \sma (G/H)_\epsilon^c \sma G/H_+ \ar[r]^-{\id \sma \id \sma \epsilon}  & 
E \sma X \sma S^V. }
\]
One consequence of Spanier--Whitehead duality is that we can factor $\tilde{\xi}$ as the 
composite
\[E \sma X \sma (G/H)_\epsilon^c \rightarrow F(G/H_+, S^V) \sma E \sma
X \rightarrow F(G/H_+, S^V \sma E \sma X).\] Since $G/H_+$ is dualizable, the second map is a 
$\pi_*$--equivalence. The first map is obtained from the duality equivalence
\[\Sigma^\infty (G/H)_\epsilon^c \rightarrow F(G/H_+, S^V)\]
by smashing with $E \sma X$ on both sides.
\vspace{-6pt}

Though $E \sma X$ is not necessarily cofibrant, $\Sigma^\infty (G/H)_\epsilon^c$ is cofibrant 
and $F(G/H_+, S^V)$ is homotopy equivalent to a cofibrant $\sWuG$--space, so the map is a 
$\pi_*$--equivalence.
\end{proof}
\vspace{-6pt}

\begin{remark}
One could also explicitly construct a homotopy inverse to the space-level map $Z(T\nu \sma X) 
\rightarrow \Map_0(G/H_+, F(S^V \sma X))$ using $V$--duality.  In our treatment, this is 
packaged up inside the machinery of \fullref{C:lyd}.
\end{remark}
\vspace{-6pt}

We will refer to $\sWuG$--spaces satisfying these equivalent conditions as ``genuinely 
equivariantly linear''. 
\vspace{-6pt}

\subsection{Refinement via the Wirthmuller isomorphism}
\vspace{-6pt}

Using the ideas that lead to the generalized Wirthmuller isomorphism, we can replace 
hypothesis (A) with a condition which does not explicitly involve $T\nu$.  To do so, we must 
first digress and discuss the passage from $\sWuG$--spaces to $\sWuH$--spaces induced by an 
inclusion $H \rightarrow G$.  There is a forgetful functor $\iota^* \co \sT_G \rightarrow 
\sT_H$.  Since $\iota^* (G/K)$ admits a triangulation as a finite $H$--CW--complex 
\cite[5.2.2]{mandell-may}, this restricts to a forgetful functor $\iota^* \co \sWuG 
\rightarrow \sWuH$.  

\begin{definition}
Given $H \monoto G$, define $\iota^* Z$ as $(\iota^*Z)(\iota^* A) =
  \iota^* (Z(A))$.
\end{definition}

Of course, not all $A \in \sWuH$ are in the image of $\iota^*$, and thus what we have really 
produced is a continuous functor from $\iota^* \sWuG$ to $H$--spaces.  To obtain an 
$\sWuH$--space, we apply the prolongation functor along the inclusion of $\iota^* \sWuG$ in 
$\sWuH$.  This process is completely analogous to the construction of the change-of-group 
functors for classical prespectra, where restriction to indexing sequences and a 
change-of-universe are necessary.  More precisely, $\iota^*$ on $\sWuG$--spaces is compatible 
under the passage to prespectra with the usual change-of-group functor there, essentially by 
construction.  That is, $\iota^* \bU Z \cong \bU \iota^* Z$.

\begin{remark}
One can show that there is a Quillen equivalence between the stable model structures on 
$\iota^* \sWuG$--spaces and $\sWuH$--spaces, by comparing each category to an appropriately 
indexed category of orthogonal spectra and using the observation that the change of universe 
functors are compatible. 
\end{remark}

Henceforth, given a $\sWuG$--space we will tacitly apply it to $H$--spaces and mean the 
corresponding $\sWuH$--space produced in the fashion above.  With this in hand, we proceed to 
revise hypothesis (A) using the Wirthmuller isomorphism.

The ``neo-classical'' construction of the Wirthmuller isomorphism given in \cite{lms} depends 
on a space-level $H$--map 
\[u \co G \sma_H X \rightarrow S^L \sma X,\]
where $G/H$ is embedded in a representation $V$, $L$ is the associated tangent
$H$--representation at the identity, and $X$ is an arbitrary $H$--space. Applying $Z$ to both 
sides, we get an $H$--map 
\[Z(G \sma_H X) \rightarrow Z(S^L \sma X).\]
Using the fact that $Z(G \sma_H X)$ is regarded as an $H$--space by forgetting down from the 
$G$--space structure and the adjunction between the forgetful functor and $\Map_H$, we get an 
induced $G$--map
\[\xi_3 \co Z(G \sma_H X) \rightarrow \Map_H(G_+, Z(S^L \sma X)).\]
Now, if we let $X$ be $S^W$ where $W \oplus L = V$ as an $H$--space, we get an $H$--map
\[Z(G \sma_H S^W) \rightarrow \Map_H(G_+,Z(S^V))\]
which corresponds to a $G$--map
\[\xi_3 \co Z(G \sma_H S^W) \rightarrow \Map_0(G/H_+, Z(S^V)).\]
Since $G \sma_H S^W$ is precisely $T\nu$ in this setting, we can
compare $\xi_3$ to $\xi_2$ (which was defined in \fullref{r:xi}).

\begin{lemma} \label{L:xi}
Under the identification of $G \sma_H S^W$ with $T\nu$, the maps $\xi_3$ and $\xi_2$ are 
$G$--homotopic.
\end{lemma}

\begin{proof}
This is essentially a consequence of the observation \cite[2.5.9]{lms} that the 
Pont\-rya\-gin--Thom map $S^V \rightarrow G \sma_H S^W$ and the Wirthmuller map $u \co G 
\sma_H S^W \rightarrow S^V$ are compatible. That is, the composite
\[S^V \rightarrow G \sma_H S^W \rightarrow S^V\]
is $H$--homotopic to the identity.  By inspection, this permits the desired comparison of 
$\xi_3$ and $\xi_2$.
\end{proof}
\eject
\begin{hyp}
A $\sWuG$--space $Z$ satisfies hypothesis (B) for the universe $U$ if the following two 
conditions hold. 
\begin{enumerate}
\item{$Z$ takes $G$--homotopy pushout squares to $G$--homotopy
    pullback squares.}
\item{Let $G/H$ embed in a representation $V$ in the universe $U$.
  Let $L$ be the tangent $H$--representation at the identity coset.
  Then for all $X \in \sWuH$, the map 
\[\xi_3 \co Z(G \sma_H X) \rightarrow \Map_H(G_+, Z(S^L \sma X))\]
is a $G$--equivalence.}
\end{enumerate} 
\end{hyp}

As one would hope, it turns out that this is equivalent to the
previous condition.  To prove this, we need a specialization of
\fullref{C:lyd}.

\begin{lemma} \label{l:lyd}
Let $Y$ and $Z$ be $W_G$--spaces such that for all $A \in \sWuG$ the prespectra $Y_U[A]$ and 
$Z_U[A]$ are $\Omega$--$G$--prespectrum.  Take $B, C \in \sWuH$.  Then given a map of 
prespectra $Y_U[S^0] \rightarrow Z_U[S^0]$ arising from a natural transformation $Y 
\rightarrow Z$ and a map of $H$--spaces $G \sma_H B \rightarrow C \sma B$, if the induced map 
of $G$--prespectra 
\[Y_U[S^0] \sma (G \sma_H B) \rightarrow F_H(G_+, Z_U[S^0] \sma C \sma
B)\] is a $\pi_*$--equivalence then the induced map of $G$--spaces
\[Y(G \sma_H B) \rightarrow \Map_H(G_+, Z(C \sma B))\]
is a weak equivalence.
\end{lemma}

\begin{proof}
The argument is similar to the third part of \fullref{C:lyd}. There is a commutative diagram 
of $G$--prespectra 
\[
\begin{CD}
Y_U[S^0] \sma (G \sma_H B) @>>> F_H(G_+, Z_U[S^0] \sma C \sma B) \\
@VVV @VVV \\
Y_U[G \sma_H B] @>>> F_H(G_+, Z_U[C \sma B])\\
\end{CD}
\]
which arises as the adjoint of the commutative diagram of $H$--prespectra:
\[
\begin{CD}
Y_U[S^0] \sma (G \sma_H B) @>>> Z_U[S^0] \sma C \sma B \\
@VVV @VVV \\
Y_U[G \sma_H B] @>>> Z_U[C \sma B]\\
\end{CD}
\]
Now the result follows from the fact that for any $H$--prespectrum $X$,
\[\Omega^\infty F_H(G_+, X) = \Map_H(G_+, \Omega^\infty X).\proved\] 
\end{proof}

\begin{thm}
A $\sWuG$--space $Z$ satisfies hypothesis (A) if and only if it satisfies hypothesis (B).
\end{thm}

\begin{proof}
We will prove that hypothesis (B) implies hypothesis (A), and that if $Z_U[A]$ is an 
$\Omega$--$G$--spectrum for all $A \in \sWuG$ then $Z$ satisfies hypothesis (B).  The fact 
that hypothesis (B) implies hypothesis (A) is an immediate consequence of \fullref{L:xi}, 
which identifies $\xi_3$ with $\xi$.

On the other hand, if $Z_U[A]$ is an $\Omega$--$G$--prespectrum for all $A \in \sWuG$ then 
hypothesis (B) holds as a consequence of the Wirthmuller isomorphism. Once again, let $E$ 
denote $Z_U[S^0]$.  Using \fullref{l:lyd}, the map
\[\xi_3 \co Z(G \sma_H X) \rightarrow \Map_H(G_+, Z(S^L \sma X))\]
is a weak equivalence if the map
\[\tilde{\xi_3} \co E \sma (G \sma_H X) \rightarrow F_H(G_+, E \sma (S^L
\sma X))\] is a $\pi_*$--equivalence.  The map $\tilde{\xi_3}$ is constructed as follows. Via 
application of $\Sigma^{\infty}$, $u$ induces a map of $H$--prespectra 
\[\mu \co G \thp_H \Sigma^{\infty} X \rightarrow \Sigma^L \Sigma^{\infty} X\] 
and this induces a map of $G$--prespectra
\[G \thp_H \Sigma^\infty X \rightarrow F_H(G_+, \Sigma^L \Sigma^\infty X)\]
which is in fact the Wirthmuller map \cite[2.6.10]{lms}. Now smashing $\mu$ by $E$ (regarded 
as an $H$--prespectrum) on both sides prior to inducing to a map of $G$--prespectra yields a 
map
\[G \thp_H (X \sma E) \rightarrow F_H(G_+, \Sigma^L (E \sma X))\]
and using the fact that $E$ is actually a $G$--prespectrum, this simplifies to
\[\tilde{\xi_3} \co E \sma (G \sma_H X) \rightarrow F_H(G_+, E \sma (S^L \sma X)).\]
The Wirthmuller isomorphism tells us this is $\pi_*$--equivalence.
\end{proof}

This provides the connection to the model-theoretic discussion.

\begin{corollary}
A $\sWuG$--space $Z$ is fibrant in the absolute stable model structure if and only if $Z$ is 
genuinely equivariantly linear (satisfies the conditions of \fullref{t:equilin}).
\end{corollary}

\subsection{Incomplete universes}\label{s:incomplete}

In the previous sections (and in the proofs of the model structures given in the appendix), 
we assume that the universe $U$ is complete. This assumption enters into our arguments when 
we employ Spanier--Whitehead duality.  In a complete universe, all finite $G$--CW--complexes 
are stably dualizable.  However, this is no longer true in incomplete universes; as we 
briefly mentioned in the introduction, the failure of Spanier--Whitehead duality in the 
trivial universe is one of the motivating factors for the use of $G$--spectra indexed on a 
complete universe.

Nonetheless, variants of our main results are valid when $U$ is not
complete.  Duality in incomplete universes has been carefully studied
by Lewis \cite{lewis}.  For our purposes, an essential result of Lewis
is that an orbit spectrum $\Sigma^\infty (G/H)_+$ is dualizable in the
stable category with respect to an incomplete universe $U$ if and only
if $G/H$ embeds in $U$.  Lewis has also carefully verified the
existence of the Wirthmuller map we use in the context of incomplete
universes.  Therefore, we can obtain variants of hypotheses (A) and
(B) which are valid in incomplete universes by restricting to orbits
$G/H$ which embed in the universe.

The other essential modification is forced by the partial failure of the equivariant version 
of Lydakis' \fullref{p:lyd}, which states that for a $W_G$--space $Z$ and any finite $G$--CW--complex $A$, there is a $\pi_*$--equivalence $Z_U[S^0] \sma 
A \htp Z_U[A]$.  The proof of this 
proposition relies on $A$ being stably dualizable, and so when working over an incomplete 
universe we must restrict to $A \in \sWuG$ which are dualizable.  We employ this proposition 
to prove both the existence of the absolute stable model structure (and the associated 
characterization of fibrant objects) as well as to show the sufficiency of our hypotheses.  
Therefore, when working over an incomplete universe we need to restrict the quantification to 
stably dualizable $A$.

\subsection{Remarks on equivariant infinite loop space theory}

These fibrancy criteria provide a conceptual understanding of the marked difference between 
equivariant infinite loop space theory for $G$ finite and for $G$ an infinite compact Lie 
group.  To be precise, we will first review the $\Gamma$--space approach to equivariant 
infinite loop space theory for a finite group $G$.  There are two obvious approaches to 
generalizing the nonequivariant theory of $\Gamma$--spaces.  A direct generalization is to 
consider $\Gamma$--$G$--spaces, which are functors from finite pointed sets to $G$--spaces.  
Alternatively, one could consider $\Gamma_G$--spaces, continuous functors from the category 
of finite pointed $G$--sets to $G$--spaces.  It turns out that these categories are 
equivalent.  This comparison was observed by Shimakawa \cite{shimakawa-note}, and is an 
example of a general fact about diagram spectra which (in the context of orthogonal 
$G$--spectra) is comprehensively discussed as part of the treatment of change-of-universe 
functors in \cite[5.1]{mandell-may}.

A $\Gamma_G$--space $X$ is ``special'' if the natural map
\[X(G/H) \rightarrow \Map_0(G/H_+, X(1))\]
is an equivariant weak equivalence.  This condition is equivalent to the usual condition on 
$\Gamma$--$G$--spaces \cite{shimakawa-note}. Associated to a $\Gamma_G$--space $X$ via 
prolongation is a $\sWuG$--space $\bP` X$.  The main theorem in this setting is that a ``very 
special'' $\Gamma_G$--space gives rise to an $\sWuG$--space which almost satisfies our 
hypothesis (A). Specifically, it only takes some homotopy pushouts squares to homotopy 
pushout squares.  However, it turns out that enough of hypothesis (A) is satisfied for the 
prespectrum $\bU \bP` X$ to be identifiable as a positive $\Omega$--$G$--prespectra.  Recall 
that a positive $\Omega$--$G$--prespectrum is a $G$--prespectrum $Y$ such that the adjoint 
structure maps $Y(V) \rightarrow \Omega^W Y(V \oplus W)$ are weak equivalences for $W$ such 
that $W^{G} \neq 0$.

\begin{remark}
Analogous to the positive stable model structure on orthogonal $G$--spectra there is a 
positive stable model structure on $\sWuG$--spaces, obtained using identical arguments to 
those presented above to construct the stable model structure.
\end{remark}

Consider hypothesis (B) in the case when $G$ is finite.  Then we
know that $\Sigma^\infty G/H_+$ is self-dual or equivalently
$\Sigma^\infty T\nu$ is the same as $\Sigma^\infty G/H_+$.  Therefore,
hypothesis (B) amounts to requiring that the map
\[Z(G \sma_H X) \rightarrow \Map_H(G/H_+, Z(X))\]
be a $G$--equivalence.  Plugging in $S^0$, we recover the ``special'' condition on the 
underlying $\Gamma_G$--space.

Now let $G$ be an infinite compact Lie group.  In this setting, we can only consider 
$\Gamma$--$G$--spaces, as finite sets do not admit interesting $G$--actions.  The dual of 
$\Sigma^\infty G/H_+$ is $G \thp S^{-L}$, and the representation sphere $S^{-L}$ is often 
nontrivial.  Even restricting to $X = S^0$ in hypothesis (B), we must consider the map
\[Z(G/H) \rightarrow \Map_H(G_+, Z(S^L)).\]
It is difficult to imagine how this equivalence could be encoded by
entirely discrete data.

Instead, these requirements strongly suggests that a reasonable domain category for the 
correct analogue of $\Gamma$--spaces must contain enough information to encode these 
dualities, and therefore most likely should contain the orbit spectra $G/H$.  Finally, it is 
worth pointing out that in general when $H$ has finite index in $G$, then $L$ is also 
trivial.  Amongst other things, this suggests that infinite loop space theory for profinite 
groups when restricting the universe to finite index subgroups should be tractable.

\appendix

\section[Model category structures on G S T]{Model category structures on $G` \sW\! \sT$}

In this section we will analyze $\sWuG$--spaces as equivariant diagram spectra.  Closely 
following \cite{mandell-may}, we will construct stable model structures on $G` \sW\! \sT$, 
compare these model structures to the stable model structure on orthogonal $G$--spectra, and 
provide a model-theoretic characterization of the fibrant objects.

\subsection[A rapid review of G-topological model categories]{A rapid review of $G$--topological model categories}

The model structures we construct on $G`\sW\!\sT$ are compatible with the enrichment in based 
$G$--spaces.  This is expressed by an appropriate variant of Quillen's SM7 axiom, as follows.  
We briefly recall from \cite[3.1.4]{mandell-may} the following definition.  Assume that we 
have a $G$--category $\sC_G$, and its associated category of $G$--maps $G`\sC$.  Let $i \co A 
\rightarrow X$ and $p \co E \rightarrow B$ be maps in $G`\sC$ and consider the map 
\[\sC_G(i^*,p_*) \co \sC_G(X,E) \rightarrow \sC_G(A,E)
\times_{\sC_G(A,B)} \sC_G(X,B)\] 
induced by $\sC_G(i, \id)$ and $\sC_G(\id,p)$ by passage to
pullbacks. 

\begin{definition}
A model category is $G$--topological if the map $\sC_G(i^*,p_*)$ is a Serre fibration (of 
$G$--spaces) when $i$ is a cofibration and $p$ is a fibration and is a weak equivalence when 
in addition either map is a weak equivalence.
\end{definition}

\subsection{The stable model structure}

The construction of the various model structures on $G`\sW\!\sT$ is mostly formal, using the 
technology developed in \cite{mandell-may} and \cite{mmss}.  In the remainder of this section 
the predominant emphasis is on recording results along with carefully verifying the specific 
variant technical lemmas necessary for this situation.  The interested reader should refer to 
the cited sources to reconstruct full arguments. 

\begin{remark}
Note that we will provide model structures only for $G`\sW\! \sT$; it isn't very useful to 
talk about such structure on $\sWuG \sT$. Nonetheless, $\sWuG \sT$ is an important device for 
encoding the compatibility of the model structure on $G`\sW\! \sT$ with the enrichment.
\end{remark}

Throughout, fix a complete universe $U$.  The first model structure we consider is the 
relative level model structure on $G`\sW\! \sT$, where by relative we mean that the 
fibrations and weak equivalences are detected only on the spheres $\{S^V\}$ for $V \in U$.

\begin{definition}
The relative level model structure on this category is defined as follows.  A map $Y \rightarrow Z$ is 
\begin{enumerate}
\item{a fibration if each $Y(S^V) \rightarrow Z(S^V)$ is an equivariant Serre fibration,}
\item{a weak equivalence if each $Y(S^V) \rightarrow Z(S^V)$ is an equivariant weak equivalence,}
\item{a cofibration if it has the left-lifting property with respect to the acyclic fibrations.} 
\end{enumerate}
\end{definition}

\begin{proposition}
The relative level model structure on $G`\sW\!\sT$ is a cofibrantly generated 
$G$--topological model structure.
\end{proposition}

\begin{proof}
The arguments are the same as \cite[3.2.4]{mandell-may}.
\end{proof}

There is an associated stable model structure.  We define $\pi_* Z$ for a $\sWuG$--space $Z$ 
by passing to the $G$--prespectrum $\bU Z$ and specifying $\pi_* Z = \pi_* \bU Z$.  Recall 
that for a subgroup $H$ of $G$ and an integer $q$, for $q \geq 0$ we define
\[\pi^H_q(\bU Z) = \colim_V \pi_q^H (\Omega^V Z(V))\]
and for $q > 0$ we define
\[\pi^H_{-q}(\bU Z) = \colim_{V \supseteq \R^q} \pi_0^H (\Omega^{V -
  \R^q} Z(V)).\]

\begin{definition}
In the stable model structure, a map is 
\begin{enumerate}
\item{a cofibration if it is a cofibration in the relative level model structure,}
\item{a weak equivalence if it is a $\pi_*$--equivalence,}
\item{a fibration if it has the right-lifting property with respect to the acyclic cofibrations (maps which are both level cofibrations and $\pi_*$--equivalences).}
\end{enumerate}
\end{definition}

We can employ the argument of \cite[3.4.2]{mandell-may} to prove that this is a model 
structure, but we need to specialize a lemma to the current situation.  Recall that $A 
\mapsto F_B A$ as a functor from $G$--spaces to $\sWuG$--spaces is defined to be left adjoint 
to the functor which is evaluation at $B$.  Concretely, we have $(F_B A) (C) = \Map_0(B,C) 
\sma A$.

There is a map \[\lambda_{V,A} \co F_{\Sigma^V A} S^V \rightarrow F_A S^0\] defined to be map such that 
\[\lambda_{V,A}^* \co \sWuG \sT (F_A S^0, X) \rightarrow \sWuG \sT (F_{\Sigma^V A} S^V, X)\]
corresponds under adjunction to 
\[X(A) \rightarrow \Omega^V X(\Sigma^V A).\]
The functors $\lambda_{V,A}$ play a key role in constructing the stable model structures, as they allow us to provide explicit descriptions of the generating cofibrations.

\begin{lemma}\label{l:lambda}
For all based $G$--CW--complexes $B$, the maps
\[\lambda_{V,A} \sma \id \co F_{\Sigma^V A} (\Sigma^V B) \cong F_{\Sigma^V A} S^V \sma B \rightarrow F_A S^0 \sma B \cong F_A B\]
are $\pi_*$--equivalences. 
\end{lemma}

\begin{proof}
We can write the specified map at a sphere $S^Z$ as
\[\Map_0(S^V \sma A, S^Z) \sma S^V \sma B \rightarrow \Map_0(A, S^Z) \sma B.\]
Rewriting, this is
\[(\Sigma^V \Omega^V \Map_0(A, S^Z)) \sma B \rightarrow \Map_0(A, S^Z) \sma B\]
and the map is the evaluation map.  First, we can assume that $B = S^0$.  It will suffice to 
show that the map is a $\pi_*$--equivalence in this case, since $B$ is a $G$--CW--complex and 
hence smashing with $B$ preserves $\pi_*$--equivalences.  But observe that $\Map_0(A,S^Z)$ is 
describing the $Z$--th space of the cotensor prespectrum $F(A,S)$, and the map in question is 
a stable equivalence because the unit $\Sigma^V \Omega^V X \rightarrow X$ is a stable 
equivalence of prespectra.
\end{proof}

\begin{remark}
Note that we could also prove this directly by induction over cell
decompositions, as is done in \cite[17.1]{mmss}, if a self-contained proof
was desired that did not require the prior work on prespectra.
\end{remark}

We need the following corollary, which trivially follows by setting $A
= S^W$ in the lemma.

\begin{corollary}\label{c:lambda}
For all based $G$--$CW$ complexes $B$, the maps
\[\lambda_{V,S^W} \sma \id \co F_{S^{V \oplus W}} \Sigma^V B \cong F_{S^{V \oplus W}} S^V \sma B \rightarrow F_{S^W} S^0 \sma B \cong F_{S^W} B\]
are $\pi_*$--isomorphisms.
\end{corollary}

We define the generating cofibrations as follows.  First, recall the sets $I$ and $J$ from \cite[3.1.1]{mandell-may}; $I$ is the set of cell cofibrations 
\[i\co (G/H \times S^{n-1})_+ \rightarrow (G/H \times D^n)_+\]
and $J$ is the set of cofibrations
\[i_0 \co (G/H \times D_n)_+ \rightarrow (G/H \times D^n \times [0,1])_+.\]
Here $H$ runs through the closed subgroups of $G$ and $n \geq 0$.

\begin{definition}
The set $FI$ is the set of all maps $F_{S^V} i$ for $i \in I$ and $V
\subset U$.  The set $FJ$ is the set of all maps $F_{S^V} j$ for $j \in J$
and $V \subset U$.
\end{definition}

We need to define the operation $f \Box g$ for maps $f$ and $g$ in order to specify the generating acyclic cofibrations.  If $i \co X \rightarrow Y$ and $j \co W \rightarrow Z$ are cofibrations, then there is a cofibration 
\[i \Box j \co (Y \sma W) \cup_{X \sma W} (X \sma Z) \rightarrow Y \sma Z.\]

\begin{definition}
Let $M\smash{\lambda_{V,S^W}}$ be the mapping cylinder of $\smash{\lambda_{V,S^W}}$.  Then 
$\smash{\lambda_{V,S^W}}$ factors as the composite of a cofibration $\smash{k_{V,W} \co 
F_{S^{V \oplus W}} S^W \rightarrow M\lambda_{V,S^W}}$ and a deformation retraction 
$\smash{r_{V,W} \co M\lambda_{V,S^W} \rightarrow F_{S^V} S^0}$.  Let $K$ be the union of $FJ$ 
and maps of the form $i \Box k_{V,W}, i \in I$.
\end{definition}

Using \fullref{c:lambda}, the arguments of \cite[3.4.2]{mandell-may} then
imply the following result.

\begin{proposition}
The stable model structure is a $G$--topological model structure on the category $G` \sW\! 
\sT$ with generating cofibrations $FI$ and generating acyclic cofibrations $K$.
\end{proposition}

\subsection[Comparison to orthogonal G-spectra]{Comparison to orthogonal $G$--spectra}

There is a pair of adjoint functors ($\bP$, $\bU$) connecting $G` \sW\! \sT$ and the category 
$G``\sI\! \sS$ of orthogonal $G$--spectra \cite{mandell-may,mmss}.  $\bU$ is the forgetful 
functor from $\sWuG$--spaces to orthogonal $G$--spectra, and $\bP$ is the prolongation 
constructed as a left Kan extension along the inclusion of domain categories 
\cite[23.1]{mmss}.  As an immediate consequence of \cite[2.14]{mandell-may}, we find that 
these functors preserve the symmetric monoidal structures.

\begin{lemma}
$\bP$ is a strong symmetric monoidal functor and $\bU$ is a lax
symmetric monoidal functor. 
\end{lemma}

Moreover, we have the expected comparison result.

\begin{thm}
The pair $(\bP,\bU)$ specifies a Quillen equivalence between the stable model category 
structure on $G` \sW\! \sT$ and the stable model category structure on $G`` \sI\! \sS$.
\end{thm}

\begin{proof}
This is virtually identical to the comparison between orthogonal $G$--spectra and 
$G$--prespectra \cite[4.16]{mandell-may}. 
\end{proof}

\subsection{The absolute stable model structure}

The level model structure used to construct the stable model structure in the previous 
section depends on evaluation at the spheres.  This makes it inconvenient to compare to 
diagram categories where the domain does not include an embedding of the spheres, for 
instance $\Gamma$--$G$--spaces.  As in the nonequivariant case, we rectify this by 
considering an ``absolute'' model structure.

\begin{definition}
The absolute level model structure on $G` \sW\! \sT$ is defined as follows.  A map $Y 
\rightarrow Z$ is
\begin{enumerate}
\item{a fibration if each $Y(A) \rightarrow Z(A)$ for $A \in \sWuG$ is
  an equivariant Serre fibration,}
\item{a weak equivalence if each $Y(A) \rightarrow Z(A)$ for $A \in
  \sWuG$ is an equivariant weak equivalence,}
\item{a cofibration if it has the left-lifting property with respect
  to the acyclic fibrations.} 
\end{enumerate}
\end{definition}

\begin{proposition}
The absolute level model structure on $G` \sW\! \sT$ is a cofibrantly generated 
$G$--topological model structure.
\end{proposition}

\begin{proof}
Again, the arguments are the same as \cite[3.2.4]{mandell-may}.
\end{proof}

There is an associated absolute stable model structure.

\begin{definition}
In the absolute stable model structure, a map is
\begin{enumerate}
\item{a cofibration if it is a cofibration in the absolute level model structure,}
\item{a weak equivalence if it is a $\pi_*$--equivalence,}
\item{a fibration if it has the right-lifting property with respect to the acyclic cofibrations.}
\end{enumerate}
\end{definition}

In order to prove that these definitions yield a model structure on the category of 
$\sWuG$--spaces, we require the full strength of \fullref{l:lambda}, which tells us that the 
maps $\lambda_{V,A}$ are $\pi_*$--equivalences.  To obtain the generating cofibrations and 
acyclic cofibrations, we enlarge $FI$ and $FJ$ by defining $F^{\prime}I$ to be the set of all 
maps $F_A i$ for $i \in I$ and $A \in \sWuG$ and $F^{\prime}J$ to be the set of all maps $F_A 
j$ for $j \in J$.  We then construct $K^{\prime}$ analogously to $K$, taking mapping 
cylinders for all maps $\lambda_{V,A}$.

\begin{proposition}
The stable model structure is a $G$--topological model structure on the category $G` \sW\! 
\sT$ with generating cofibrations $F^{\prime}I$ and generating acyclic cofibrations 
$K^{\prime}$.
\end{proposition}

\begin{proof}
This follows the proof of \cite[17.2]{mmss}, modified slightly in light of the proof of \cite[3.4.2]{mandell-may}.
\end{proof}

In the course of this proof we obtain the following analogue of \cite[3.4.8]{mandell-may}.

\begin{proposition}
A $\sWuG$--space $Z$ is fibrant in the absolute stable model structure if and only if for all 
$A \in \sWuG$ and $W \in U$ the structure map 
\[Z(A) \rightarrow \Omega^W Z(S^W \sma A)\] 
is a weak equivalence.
\end{proposition}

Finally, we can compare the absolute and relative stable model
structures.  It is clear that the identity functor is the left adjoint
in a Quillen pair relating the two model structures. 

\begin{thm}
The identify functor is the left adjoint of a Quillen equivalence between the category of 
$\sWuG$--spaces with the relative stable model structure and the category of $\sWuG$--spaces 
with the absolute stable model structure. 
\end{thm}

\subsection{Ring and module spectra}

We can leverage the results of \cite{mandell-may,mmss} to lift the stable model structure to 
categories of ring and module $\sWuG$--spaces.  The key ingredient in these lifting results 
is the verification of the monoid axiom and the pushout-product axiom.

To verify these, we need the following technical lemma.

\begin{lemma}\label{l:techlem}
Let $Y$ be a $\sWuG$--space such that $\pi_*(Y) = 0$.  Then $\pi_*(F_V S^V \sma Y) = 0$ for 
any $V$. 
\end{lemma}

\begin{proof}
The conclusion follows immediately from the counterpart for orthogonal $G$--spectra
\cite[3.7.2]{mandell-may} upon applying the prolongation functor $\bP$ to $\sWuG$--spaces, 
just as in \cite[12.3]{mmss}. 
\end{proof}

Now, the same chain of arguments given in \cite[3.7]{mandell-may}
allows us to verify the following two results. 

\begin{proposition}[Monoid axiom]
For any acyclic cofibration $i \co X \rightarrow Y$ of $\sWuG$--spaces and any $\sWuG$--space 
$Z$, the map $i \sma \id \co X \sma Z \rightarrow Y \sma Z$ is a $\pi_*$--isomorphism and an 
$h$--cofibration.  This holds for cobase changes and sequential colimits of such maps as 
well. 
\end{proposition}

Note that by $h$--cofibration we mean a map satisfying the homotopy extension property as 
opposed to a model theoretic cofibration.

\begin{proposition}[Pushout-product axiom]
If $i \co X \rightarrow Y$ and $j \co W \rightarrow Z$ are cofibrations of $\sWuG$--spaces 
and $i$ is a $\pi_*$--isomorphism, then the cofibration $$i \Box j \co (Y \sma W) \cup_{X \sma 
W} (X \sma Z) \rightarrow Y \sma Z$$ is a $\pi_*$--isomorphism. 
\end{proposition}

As an immediate consequence, we have the following version of
\cite[3.7.6]{mandell-may}. 

\begin{thm}
Let $R$ be a ring $\sWuG$--space.
\begin{enumerate}
\item{The category of $R$--module $\sWuG$--spaces is a cofibrantly generated proper
  $G$--topological model category, with weak equivalences and
  fibrations created in the stable model structure on the category of $\sWuG$--spaces.}
\item{If $R$ is commutative, the category of $R$--algebra
  $\sWuG$--spaces is a cofibrantly generated right proper
  $G$--topological model category with weak equivalences and fibrations
  created in the stable model structure on the category of $\sWuG$--spaces.}
\end{enumerate}
\end{thm}

\begin{remark}
Following the outline above, one can also obtain a version of this
theorem by lifting the absolute stable model structure.  This variant
of the theorem requires a slightly stronger version of \fullref{l:techlem}, obtained by an equivariant version of the argument for
\cite[17.6]{mmss}. 
\end{remark}

\section[The failure of the approximation theorem for G = S1]{The failure of the approximation theorem for $G = S^1$}

An equivariant approximation theorem would purport to show that an
appropriate map $C(V) \rightarrow \Omega^V S^V$ was a group
completion, where $C(V)$ is the configuration space of points of $V$.
Note that this formulation would actually be correct only for $V$ such
that $V^G$ is nonzero.  The statement is somewhat more complicated when
$V^G = 0$, as then there is no addition. 

This counterexample is due to Segal \cite{equi-segal}.  Let $G = S^1$ and let $V = \R^3$ 
where $G$ acts by rotation around the $z$--axis. The inclusion of the axis gives a 
cofibration $S^1 \rightarrow S^V$ of $G$--spaces, and the cofiber $S^V / S^1$ is 
$G$--homeomorphic to $\Sigma^2 G_+$.  To see this, let $G$ be the unit circle in the 
$(x,y)$--plane with the disjoint basepoint at the origin, and parameterize $\R^3$ as $S^1 
\times [0,\infty) \times \R / (t,0,s) \sim (t^{\prime}, 0, s)$. 

Choose a particular $G$--space $X$.  There results a $G$--fibration sequence
\[\Map_0(\Sigma^2 G_+, X) \rightarrow \Map_0(S^V,X) \rightarrow
\Map_0(S^1,X).\] Passing to fixed-points, this remains a fibration sequence.  Evaluating the 
fixed-points of the terms in the sequence, we find \begin{gather*}\Map_0(S^1,X)^G \htp 
\Omega(X^G)\\
 \Map_0(\Sigma^2 G_+, X)^G
\cong \Map_0(G_+, \Omega^2 X)^G \cong \Omega^2 X.\tag*{\hbox{and}}\end{gather*} Thus, we have 
the fibration sequence
\[\Omega^2 X \rightarrow (\Omega^V X)^G \rightarrow \Omega(X^G).\]
Finally, take $X = S^V$.  Then $X^G = S^1$ and $\Omega S^1 \htp \Z$. The sequence splits and 
so we have $(\Omega^V S^V)^G \htp \Z \times \Omega^2 S^3$.  But $G$--fixed points of the 
configuration space $C(V)^G$ are the same as the configuration space $C(\R^1)$, and the usual 
group completion theorem tells us that the group completion of this is $\Omega S^1$.  
Therefore, we have a contradiction.

\section[The trouble with Gamma--S1--spaces]{The trouble with $\Gamma$--$S^1$--spaces}

In naive analogy with the situation for $G$ finite, one might hope that there is some 
condition on a $\Gamma$--$S^1$--space $F$ which would guarantee that the $W_{S^1}$--space 
$\bP F$ obtained by prolongation would be a positive $\Omega$--$S^1$--spectrum 
\cite{equi-segal,shimakawa}.  However, we will show that no such condition can exist by 
studying the counit of the $(\bP, \bU)$ adjunction.

Any satisfactory condition would certainly be satisfied by a $\Gamma$--$S^1$--space obtained 
by forgetting from a genuinely equivariantly linear $W_{S^1}$--space.  Consider the case in 
which we begin with a fibrant $W_{S^1}$--space $X$ such that the prespectrum $\bU X$ is 
connective.  Let $Y$ be the $\Gamma$--$S^1$--space obtained from $X$ via the forgetful 
functor, and denote by $\tilde{X}$ the prolongation $\bP Y$.  The counit of the adjunction 
gives us a map $\tilde{X} \rightarrow X$.  We will compare $\tilde{X}$ and $X$, and show that 
in fact they will almost never be stably equivalent.  As a consequence, there can be no 
condition on a $\Gamma$--$S^1$--space which will guarantee that its prolongation is a 
positive $\Omega$--$S^1$--spectrum.

We will proceed by comparing the associated $G$--prespectra $\bU \tilde{X}$ and $\bU X$.  
Abusing notation, we will also refer to these $G$--prespectra as $\tilde{X}$ and $X$.  First, 
let $H \subset S^1$ be a finite subgroup.  $Y$ determines a $\Gamma$--$H$--space $Y_H$ via 
the forgetful functor, and there is an associated $\Gamma_H$--space which we will also denote 
$Y_H$.  Now, know that $X$ satisfies hypothesis (B) and therefore $Y_H$ is very special.  
Therefore, $\bP Y_H$ is fibrant and there is an equivalence $\tilde{X}(S^V)^H \htp X(S^V)^H$.

\begin{remark}
The previous observation is the starting point for a comparison of cyclic $\Gamma$--spaces 
and the $\sF$--model structure on connective $S^1$--spectra.  Here the weak equivalences on 
$S^1$--spectra are taken to be the $\sF$--equivalences \cite[4.6.5]{mandell-may}, where $\sF$ 
denotes the family of finite subgroups of $S^1$.  We intend to discuss this comparison 
elsewhere.
\end{remark}

Now consider the $S^1$ fixed-points of $\tilde{X}(S^V)$.

\begin{lemma}
There is a weak equivalence $\tilde{X}(S^V)^{S^1} \htp X(S^{V^{S^1}})^{S^1}$.
\end{lemma}

\begin{proof}
By definition, $\tilde{X}(S^V)$ is the coend 
\[\int_{\Gamma} X(n) \times (S^V)^n.\]
Since $S^1$ is infinite, observe that $\tilde{X}(S^V)^{S^1}$ is in fact the same as 
\[\int_{\Gamma} X(n)^{S^1} \times (S^{V^{S^1}})^n.\]
This implies that $\tilde{X}(S^V)^{S^1} \simeq X(S^{V^{S^1}})^{S^1}$.
\end{proof}

As an consequence, observe that there is a levelwise weak equivalence of prespectra 
$\tilde{X}^{S^1} \rightarrow X^{S^1}$.  This observation allows us to obtain a precise 
description of the $G$--prespectrum $\tilde{X}$.  Recall that $E \sF$ is the classifying 
space of the family of finite subgroups of $S^1$, so that $(E\sF_+)^H \htp S^0$ and 
$(E\sF_+)^{S^1}$ contractible.  

\begin{proposition}
There is a zig-zag of levelwise weak equivalences of $G$--prespectra between $\tilde{X}$ and 
the following homotopy pushout:
\[
\begin{CD}
\iota_* X^{S^1} \sma E\sF_+ @>>> X \sma E\sF_+ \\
@VVV @VVV \\
\iota_* X^{S^1} @>>> \bar{X} \\
\end{CD}
\]
\end{proposition}

\begin{proof}
Since $\iota_* X^{S^1} (V) \cong X(V^{S^1}) \sma S^{V-V^{S^1}}$, the structure maps of $X$ 
induce the map in the top righthand corner.  To compute the homotopy pushout, since the level 
model structure on $\sP$ is left proper we can take the actual pushout in the diagram 
obtained by replacing $\iota_* \smash{X^{S^1}} \sma E\sF_+$ by a cofibrant $G$--prespectrum 
and the map $\iota_* \smash{X^{S^1}} \sma E\sF_+ \rightarrow X \sma E\sF_+$ by a cofibration.  
The fixed-point functor commutes with this pushout, since one leg of the diagram is a 
cofibration and hence a closed inclusion \cite[3.1.6]{mandell-may}.  When we apply 
$\smash{(-)^{S^1}}$, the top row is contractible and so there is an equivalence $(\iota_* 
\smash{X^{S^1})^{S^1} \simeq \bar{X}^{S^1}}$.  When we apply $(-)^{H}$ for $H$ a finite 
subgroup of $S^1$, the left column becomes an equivalence and so we have $X^H \simeq 
\bar{X}^H$.  To obtain the connection between $\tilde{X}$ to $\bar{X}$, we consider the 
analogous homotopy pushout:
\[
\begin{CD}
\iota_* \tilde{X}^{S^1} \sma E\sF_+ @>>> \tilde{X} \sma E\sF_+ \\
@VVV @VVV \\
\iota_* \tilde{X}^{S^1} @>>> \tilde{X}^{\prime} \\
\end{CD}
\]
The map $\tilde{X} \rightarrow X$ induces a map of homotopy pushouts, which is a weak 
equivalence at each corner by previous discussion. Therefore there is an equivalence 
$\tilde{X}^{\prime} \rightarrow \bar{X}$.  There is also a map from $\tilde{X}^{\prime} 
\rightarrow \tilde{X}$ obtained from the natural maps $\iota_* \smash{\tilde{X}^{S^1}} 
\rightarrow \tilde{X}$ and $\tilde{X} \sma E\sF_+ \rightarrow \tilde{X}$.  This map clearly 
becomes an equivalence upon application of $(-)^H$ for $H \subset S^1$ finite.  Applying 
$\smash{(-)^{S^1}}$ to the associated pushout diagram, the top row is contractible and on the 
bottom we obtain the composite 
\[(\iota_* \tilde{X}^{S^1})^{S^1} \rightarrow
(\tilde{X}^{\prime})^{S^1} \rightarrow \tilde{X}^{S^1},\] which is a weak equivalence.  Since 
$(\iota_* \tilde{X}^{S^1}) \rightarrow (\tilde{X}^{\prime})^{S^1}$ is also a weak 
equivalence, the map $\smash{(\tilde{X}^{\prime})^{S^1} \rightarrow \tilde{X}^{S^1}}$ must be 
a weak equivalence.
\end{proof}

Applying the spectrification functor $L$ we can use this description to compute the 
$S^1$--fixed points of $L\tilde{X}$.  Recall that $L$ is the left adjoint in a Quillen 
equivalence between $\sP$ with the stable model structure and $\sS$ with the generalized 
cellular model structure \cite[4.2.9]{mandell-may}.  

\begin{lemma}
The spectrum $(L\tilde{X})^{S^1}$ is weakly equivalent to 
$(L(\iota_* X^{S^1}))^{S^1}$. 
\end{lemma}

\begin{proof}
The pushout square describing $\bar{X}$ as a $G$--prespectrum is taken to a pushout of 
$G$--spectra by $L$.  For this calculation, it is convenient to assume that we have replaced 
both maps in the original square by cofibrations when computing $\bar{X}$.  Then since the 
maps in the pushout are levelwise cofibrations and hence stable cofibrations, $L$ takes them 
to cofibrations of $G$--spectra.  This implies that the resulting square of $G$--spectra is 
in fact a homotopy pushout.  Therefore it is also a homotopy pullback square.  The homotopy 
pullback can be computed by taking the actual pullback in the square  
\[
\begin{CD}
Q @>>> Y \\
@VVV @VVV \\
Y^{\prime} @>>> L\bar{X} \\
\end{CD}
\]
where $Y \rightarrow L\bar{X}$ and $Y^{\prime} \rightarrow L\bar{X}$ are fibrations, $Y \htp 
L(X \sma E\sF)$, $Y^{\prime} \htp L(\iota_* \smash{X^{S^1}})$, and $Q \htp \smash{L((\iota_* 
X)^{S^1}} \sma E\sF)$.  Applying $\smash{(-)^{S^1}}$, we observe that $\smash{Q^{S^1}}$ and 
$\smash{Y^{S^1}}$ are contractible. Therefore, we obtain an equivalence $\smash{(L(\iota_* 
X^{S^1}))^{S^1} \htp (L\bar{X})^{S^1}}$, and this implies the result.
\end{proof}

However, it is rarely the case that such an equivalence holds for an arbitrary 
$S^1$--spectrum $X$.  For instance, the equivariant tom-Dieck splitting \cite[5.11.1]{lms} 
tells us that such an equivalence does not hold for suspension spectra.

\begin{remark}
In the unstable setting, it is possible to obtain a model category structure which captures 
equivariant $S^1$--homotopy theory by gluing together a simplicial set and a cyclic set 
\cite{blumberg}.  The previous discussion can be interpreted as a demonstration that the 
stable analogue of this gluing argument fails.  The natural approach would be to attempt to 
decompose an $S^1$--spectrum $X$ into the nonequivariant $S^1$--fixed points 
$\smash{X^{S^1}}$ and the ``cyclic'' part $X \sma E\sF_+$.  To recover $X$, one would then 
glue $X \sma E\sF_+$ to the $S^1$--spectrum induced from $\smash{X^{S^1}}$, $\iota_* 
\smash{X^{S^1}}$.  But the argument above shows that we cannot recapture the $S^1$--fixed 
points in this fashion.  This is perhaps not surprising in light of the significant 
difference in complexity between the diagrams representing unstable equivariant spaces 
\cite{elmendorf} and the diagrams representing stable equivariant objects 
\cite{schwede-shipley}.
\end{remark}

\bibliographystyle{gtart}
\bibliography{link}

\end{document}